\documentclass[11pt]{amsart}
\usepackage{tikz-cd}
\usepackage{enumitem}
\setlist[description]{leftmargin=\parindent,labelindent=\parindent}
\usetikzlibrary{matrix,arrows,decorations.pathmorphing}
\usepackage{url}
\usepackage{cite}
\usepackage{fullpage}
\usepackage{hyperref}
\usepackage{verbatim}
\usepackage{comment}
\usepackage{fancyvrb}
\usepackage{fvextra}
\usepackage{caption}
\usepackage{stmaryrd}
\usepackage{appendix}

\allowdisplaybreaks

\usepackage{amsmath}         
\usepackage{mathtools}       
\usepackage{amssymb}   
\usepackage{amsthm}
\usepackage{mathrsfs} 

\usepackage{tikz-cd}
\usepackage{bbm}

\usepackage{hyperref}
\hypersetup{
	colorlinks=true,
	linkcolor=black,
	citecolor=blue,
	urlcolor=blue}

\newcommand{\Z}{\mathbb{Z}}
\newcommand{\Q}{\mathbb{Q}}

\newcommand{\C}{\mathbb{C}}

\newcommand{\Ac}{\mathcal A}
\newcommand{\Bc}{\mathcal B}

\newcommand{\Dc}{\mathcal D}
\newcommand{\Ec}{\mathcal E}

\newcommand{\Grm}{\mathrm G}
\newcommand{\Hc}{\mathcal H}

\newcommand{\Jc}{\mathcal J}

\renewcommand{\Mc}{\mathcal M}

\newcommand{\Qc}{\mathcal Q}
\newcommand{\Rc}{\mathcal R}

\newcommand{\Wc}{\mathcal W}
\newcommand{\Xc}{\mathcal X}

\newcommand{\Zc}{\mathcal Z}

\newcommand{\Tor}{\operatorname{Tor}}

\newcommand{\vol}{\operatorname{vol}}

\newcommand{\End}{\operatorname{End}}

\newcommand{\Sp}{\operatorname{Sp}}

\newcommand{\CH}{\mathsf{CH}}
\newcommand{\Jac}{\operatorname{Jac}}
\newcommand{\NL}{\operatorname{NL}}

\newcommand{\taut}{\operatorname{taut}}

\title[]{Noether-Lefschetz cycles on the moduli space of abelian varieties}
\author{Aitor Iribar L\'opez}
\date{}

\newtheorem{thm}{Theorem}
\newtheorem{prop}[thm]{Proposition}
\newtheorem{lem}[thm]{Lemma}
\newtheorem{cor}[thm]{Corollary}

\newtheorem{conj}{Conjecture}

\theoremstyle{definition}
\newtheorem{defn}{Definition}

\theoremstyle{remark}

\newtheorem{rmk}[thm]{Remark}
\newtheorem{question}{Question}

\usepackage{wrapfig}

\makeatletter
\def\@defaultbiblabelstyle#1{[#1]}
\makeatother

\begin{document}
\baselineskip=17pt
\footskip=1.5\normalbaselineskip

\begin{abstract}
    The locus of non-simple abelian varieties in the moduli space of principally polarized abelian varieties gives rise to Noether-Lefschetz cycles. We study their intersection theoretic properties using the tautological projection constructed in \cite{CMOP24}, and show that projection defines a homomorphism when restricted to cycles supported on that locus. Using Hecke correspondences and the pullback by Torelli we prove that $[\mathcal {A}_1 \times \mathcal A_{g-1}]$ is not tautological in the sense of \cite{vdG99} for $g=12$ and $g\geq 16$ even. We also explore the connections between Noether-Lefschetz cycles and the Gromov-Witten theory of a moving elliptic curve.
\end{abstract}

\maketitle

\tableofcontents

\section{Introduction}

\subsection{Moduli of abelian varieties}

Given an abelian variety $X$ of dimension $g$ and an ample line bundle $L$, there is a polarization
$$
\begin{alignedat}{2}
    \theta_L : X \, && \longrightarrow &\,  X^\vee \\
    x \, && \mapsto &\,  t_x^*L \otimes L^{-1}
\end{alignedat}
$$
with kernel of the form
$$
 \ker (\theta_L) = (\Z/d_1 \times \ldots \times \Z/d_g)^2
$$
for a unique list of positive integers $\delta = (d_1, \ldots, d_g)$ satisfying $d_i \mid d_{i+1}$ so that $\left(\prod d_i\right)^2 = \deg (\theta_L)$.
The list $\delta$ is the \emph{polarization type}. For fixed $g$ and $\delta$, the collection of all such polarized abelian varieties
\begin{equation*}\label{eqn: definition Agd}
    \Ac_{g, \delta} = \left\lbrace \begin{array}{c}
         (X,\theta) : X\text{ is an abelian variety of dimension }g\\
         \text{ and }\theta:X \longrightarrow X^\vee\text{ is a polarization of type }\delta
    \end{array}\right\rbrace
\end{equation*}
has the structure of a smooth Deligne--Mumford stack of dimension $\binom{g+1}{ 2}$, see \cite{Mum71}. We refer to \cite{CF91} and \cite{BL04} for foundational aspects of these spaces.


Over the complex numbers, we can identify $\Ac_{g,\delta}$ with the stack quotient of the Siegel upper half space $\Hc_g$ by an arithmetic subgroup of the symplectic group $\Sp(2g, \Q)$. As a moduli stack of abelian schemes, there exists a universal polarized abelian scheme with a zero section:
$$
\begin{tikzcd}
	{\Xc_{g, \delta}} \\
	{\Ac_{g, \delta}}
	\arrow["\pi"', from=1-1, to=2-1]
	\arrow["s"', bend right, from=2-1, to=1-1]
\end{tikzcd}\, .
$$
The pullback under $s$ of the relative sheaf of differentials $\Omega_\pi$ is the \emph{Hodge bundle}, which we denote by $\mathbb E_{g, \delta}$, or by $\mathbb E_g$. Alternatively, $\mathbb E_g = \pi_*\Omega_\pi$.

The polarization $\theta$ is \emph{principal} if $\delta=(1,\ldots ,1)$; we say that $\theta$ is a \emph{principal polarization}. The moduli space of principally polarized abelian varieties $\Ac_g$ will be our main focus, but the moduli spaces of higher degree polarizations will play an important role.

\subsection{Noether-Lefschetz loci in $\Ac_g$}

The N\'eron-Severi group of an abelian variety is a free abelian group, whose rank is the \emph{Picard number} of $X$. A general point in $\Ac_{g,\delta}$ has Picard number $1$, and the locus where the Picard number is at least $2$ is known as the \emph{Noether-Lefschetz locus}. The Noether-Lefschetz locus is a countable union of closed and reduced subvarieties of $\Ac_{g,\delta}$, whose general element is either
\begin{itemize}
    \item an abelian variety with a non-trivial abelian subvariety, or
    \item an abelian variety with real multiplication.
\end{itemize}
see \cite{DL90}. We will focus on the first case.


Given an integer $u \leq g/2$ and a polarization type $\delta = (d_1, \ldots , d_u)$ we consider the Noether-Lefschetz subvarieties of $\Ac_g$
$$
\NL_{g,\delta} = \left \lbrace (A, \theta) \in \Ac_g \left| \begin{array}{c}
    A\text{ has an abelian subvariety }B \\
    \text{and }\theta|_{ B}\text{ is of type }\delta 
\end{array}\right.\right\rbrace.
$$
An abelian variety is \emph{simple} if it does not have any non-trivial abelian subvariaty.

If $\delta = (d)$ (so $u=1$), we will use the notation $\NL_{g,d}$ instead.

We show in Section \ref{section: 22} that these are the irreducible components of the Noether-Lefschetz locus given by non-simple abelian varieties. For example, if $\delta = (\underbrace{1,\ldots ,1}_{u\text{ times}})$ then $\NL_{g, \delta}$ is the image of the natural product map
\begin{equation}\label{eqn: product maps}
    \Ac_u \times \Ac_{g-u} \longrightarrow\Ac_g.
\end{equation}

More generally, the locus $\NL_{g,\delta}$ will be the image of a twisted product map, as we will explain in Section \ref{section: 22}.

\subsection{The tautological rings}

The Chern classes of the Hodge bundle
$$
\lambda_i = c_i(\mathbb E_{g,\delta})
$$
generate the Chow ring of $\Ac_{g,\delta}$ for $g \leq 3$ and $\delta = (1,\ldots , 1)$ (see \cite{vdG98}). In general, the subring
$$
\mathsf{R}^*(\Ac_{g,\delta}) : = \operatorname{im} \left(\Q[\lambda_1, \ldots , \lambda_g] \longrightarrow \CH^*(\Ac_{g,\delta})\right)
$$
generated by the $\lambda_i$ will be a strict subset of the Chow ring called the \emph{tautological ring} of $\Ac_{g,\delta}$. The following result follows from \cite[Theorem 3]{vdG99}:
\begin{thm}\label{thm: tautological rings}
    Let $\delta=(d_1, \ldots , d_g)$ be a polarization type. Then
    $$
    \mathsf{R}^*(\Ac_{g,\delta}) = \frac{\Q[\lambda_1, \ldots , \lambda_g]}{\langle c(\mathbb E_{g,\delta})c(\mathbb E_{g,\delta}^\vee)-1, \, \lambda_g \rangle}\, .
    $$
\end{thm}
The proof appears in Section \ref{section: 21}. The tautological ring is in fact isomorphic to the cohomology ring of the Grassmannian of Lagrangian hyperplanes in a vector space of dimension $2(g-1)$, and therefore it is Gorenstein with socle of degree $\binom{g}{ 2}$ and independent of $\delta$, see \cite{vdG99}.

\subsection{The action of the Hecke algebra}
In \cite{COP}, it was proven that the class $[\Ac_1 \times \Ac_{5}]\in \CH^5(\Ac_6)$ is not tautological. Their method consists of computing the pullback of $[\Ac_1 \times \Ac_{5}]$ under the Torelli morphism
$$
\Tor : \Mc_{g}^{ct} \longrightarrow \Ac_g\,,
$$
that sends a curve to the group of multidegree $0$ line bundles on it, where $\Mc_{g}^{ct}$ is the moduli space of stable curves of \emph{compact type}. Then, they use the completeness of Pixton's relations \cite{PPZ15} on $\overline{\Mc}_{6}$, proven by Canning, Schmitt and Larson \cite{CLS} and computer computations to check that $\Tor^*([\Ac_1 \times \Ac_{5}])$ does not agree with the pullback of its tautological projection \eqref{eqn: taut of A1 Ag-1}, which we define in the next section. We extend their results using a different strategy, which we explain now.

For each $n \geq 1$, consider the stack
$$
\mathcal B_{g,n} = \left\lbrace  f: (X,\theta_X) \longrightarrow(Y, \theta_Y) \left| \begin{array}{c}
    (X, \theta_X), (Y, \theta_Y) \in \Ac_g, f\text{ is an} \\
    \text{homomorphism, and }f^*\theta_Y = n\theta_X 
\end{array}\right.\right\rbrace.
$$
There is a correspondence
$$
\begin{tikzcd}
	& \Bc_{g,n} \\
	{\Ac_g} && {\Ac_g}
	\arrow[from=1-2, to=2-1]
	\arrow[from=1-2, to=2-3]
\end{tikzcd}
$$
where the morphisms are just the forgetful maps. By pullback from the left and pushforward to the right, they give rise to Hecke operators $T_n : \CH^{g-1}(\Ac_g) \longrightarrow\CH^{g-1}(\Ac_g)$. In Section \ref{section: 4} we show that the Hecke operators preserve tautological classes but shift the polarization type of the $\NL_{g,d}$. The pullback of $\NL_{g, 2}$ under the Torelli morphism, restricted to the locus of smooth curves
\begin{equation}\label{eqn: Torelli morphism}
\Tor : \mathcal M_{g} \longrightarrow \Ac_g
\end{equation}
is the bielliptic locus\footnote{The bielliptic locus has already been used in \cite{GP03, vZ16} to show the existence of non-tautological cycles on moduli spaces of curves}, which is shown to be non-tautological for $g =12$ in \cite{vZ16}, and $g \geq 16$ even in \cite{ACHLMT24}. This proves the following:
\begin{thm}\label{thm: not tautological}
    For $g=12$ and $g \geq 16$ even, we have that $[\Ac_1 \times \Ac_{g-1}]\not \in \mathsf{R}^{g-1}(\Ac_g)$.
\end{thm}
Note that we only need to work on the moduli space of smooth curves, where the pullback under Torelli is easier to compute.

\subsection{The tautological projection}

Canning, Molcho, Oprea and Pandharipande have constructed in \cite{CMOP24} a way to associate to any cycle class a tautological class in a natural way, which we call its \emph{tautological projection}. They calculate the tautological projection of product cycles $[\Ac_{g_1} \times \ldots \times \Ac_{g_l}]$. Finding formulas for the tautological projection of other natural families of cycles is an interesting question.

We explain their construction briefly. Let $\overline{\Ac}_g$ be a smooth toroidal compactification of $\Ac_g$, as constructed in \cite{AMRT10}, and define the $\lambda_g$-pairing on $\Ac_g$:
$$
\CH^k(\Ac_g) \times \CH^{\binom{g}{2}-k}(\Ac_g) \longrightarrow \CH^{\binom{g}{2}}(\Ac_g) \longrightarrow \Q
$$
by the formula $(\alpha, \beta) \mapsto \int_{\overline{\Ac}_g} \bar{\alpha}\cdot \bar{\beta}\lambda_g$, where $\bar{\alpha}$ and $\bar{\beta}$ are extensions of the cycles $\alpha$ and $\beta$. The pairing is well defined because
$$
\lambda_g |_{\overline{\Ac}_g \smallsetminus \Ac_g} = 0\, ,
$$
see \cite[Theorem 3]{CMOP24}.

Unlike the case of curves, the restriction of the $\lambda_g$ pairing to the tautological ring of $\Ac_g$ is perfect; this was shown in \cite{vdG99}.
\begin{defn}[\protect{\cite[Definition 4]{CMOP24}}]
    The \emph{tautological projection} $\taut : \CH^k(\Ac_g) \longrightarrow\mathsf{R}^k(\Ac_g)$ is the composition:
    $$
        \CH^k(\Ac_g) \longrightarrow (\CH^{{\binom{g}{2}}-k} (\Ac_g))^\vee \longrightarrow (\mathsf{R}^{{\binom{g}{2}}-k}(\Ac_g))^\vee \stackrel{\cong}{\longrightarrow} \mathsf{R}^k(\Ac_g)\, ,
    $$
    where the first and third arrows are given by the $\lambda_g$-pairing, and the second one is the dual to the inclusion.
\end{defn}
Unraveling the definitions, we see that $\taut(\alpha)$ is the unique tautological class such that
    $$
    \int_{\overline{\Ac}_g} \alpha \cdot \Lambda \cdot \lambda_g = \int_{\overline{\Ac}_g} \taut(\alpha) \cdot \Lambda \cdot \lambda_g
    $$
for any tautological class $\Lambda$.

The tautological projections of the product cycles $[\Ac_u \times \Ac_{g-u}]$ were calculated in \cite[Theorem 6]{CMOP24}. For example,
\begin{align}\label{eqn: taut of A1 Ag-1}
    \taut([\Ac_1 \times \Ac_{g-1}]) &= \frac{g}{6|B_{2g}|} \lambda_{g-1}\, ,\\
    \taut([\Ac_2 \times \Ac_{g-2}]) &= \frac{g(g-1)}{360 |B_{2g}||B_{2g-2}|} \lambda_{g-1} \lambda_{g-3}\, ,
\end{align}
where $B_i$ is the $i$-th Bernoulli number. In Section \ref{section: 31}, we extend their result to the Noether-Lefschetz locus given by non-simple abelian varieties:
\begin{thm}\label{thm: projection of NLD}
    If $\delta = (d_1, \ldots , d_u)$ is a polarization type and $d = d_1 \cdot \ldots \cdot d_u$ then
    $$
    \taut([\NL_{g,\delta}]) = c \cdot \taut([\Ac_u \times \Ac_{g-u}])\, ,
    $$
    where
    $$
    c = d_1^{2u-2}d_2^{2u-6}\ldots d_u^{-2u+2}\prod_{1\leq i<j\leq u} \prod_{p \mid d_j/d_i}\frac{(1-p^{-2(j-i)})}{(1-p^{-2(j-i+1)})}d^{2(g-u)+1}\prod_{j=1}^u\prod_{p \mid d_j}(1-p^{-2(j +g-2u))})\, ,
    $$
    and the products are over primes.
\end{thm}

In particular, we have\footnote{when $g=2$, $\NL_{g,d}$ is the Humbert surface of discriminant $d^2$, and its class in $\CH^1(\Ac_2)$ was computed in \cite[Theorem 8.1.]{vdG82}}
$$
    \taut([\NL_{g,d}])=\frac{gd^{2g-1}}{6 |B_{2g}|} \cdot \prod_{p \mid d}(1-p^{2-2g}) \lambda_{g-1}\, ,
$$
and $\taut([\NL_{g, (d_1, d_2)}])$ is
$$
    \left(\frac{g(g-1)d_1^{2g-1}d_2^{2g-5}}{360|B_{2g}||B_{2g-2}|} \prod_{p \mid d_1}(1-p^{6-2g}) \prod_{p \mid d_2}(1-p^{4-2g})\prod_{p \mid d_2/d_1}\frac{(1-p^{-2})}{(1-p^{-4})}\right)\lambda_{g-1}\lambda_{g-3}\, .
$$

Greer and Lian \cite{GL24, GLS24} consider the stacks
$$
\widetilde{\NL}_{g,d} = \left\lbrace f:E \longrightarrow(A,\theta)\left|\begin{array}{c}(A, \theta) \in \Ac_g, E\text{ is an elliptic curve and }\\f\text{ is an homomorphism such that }\deg(f^*\theta)=d\end{array}\right.\right\rbrace\, ,
$$
which have several connected components that correspond to our $\NL_{g,e}$ for $e \mid d$. They conjecture that the pushforward to $\Ac_g$ of their fundamental classes should be the Fourier coefficients of a modular form of weight $2g$. As a consequence of Theorem \ref{thm: projection of NLD}, we can partially verify the conjecture:
\begin{cor}\label{cor: modular form}
Define $[\widetilde{\NL}_{g,0}]$ as $ \frac{(-1)^g}{24}\lambda_{g-1}$. Under the change $q=e^{2\pi i \tau}$, we have that
$$
\taut\left( \sum_{d \geq 0}[\widetilde{\NL}_{g,d}]\, q^d\right) = \frac{(-1)^g}{24}E_{2g}(\tau)\lambda_{g-1}\, ,
$$
where $E_{2g}$ is the Eisenstein modular function, of weight $2g$.
\end{cor}

In Section \ref{section: 4} we introduce the $\mathbb Q$-vector space $\mathbf{NL}_g \subset\CH^{g-1}(\mathcal A_g)$ generated by the cycles $[\widetilde{\operatorname{NL}}_{g,d}]$ for $d \geq 0$. A. Pixton has conjectured that $\mathbf{NL}_g$ has the same dimension as the space $\operatorname{Mod}_{2g}$ of modular forms of weight $2g$ for the group $\operatorname{SL}_2(\mathbb Z)$. In light of Corollary \ref{cor: modular form}, we propose the following strengthening of Pixton's conjecture.

\begin{conj}\label{conj: modularity}
    There exists a perfect pairing
    $$
    \Phi : \operatorname{Mod}_{2g} \times \mathbf{NL}_g \longrightarrow \mathbb Q
    $$
    such that, for any modular form $f$ of weight $2g$,
    $$
    \sum_{d \geq 0} \Phi(f, [\widetilde{\NL}_{g,d}]) q^d
    $$
    is the Fourier expansion of $f$ at $i \infty$, where $q= e^{2\pi i \tau}$.
\end{conj}
By \cite[Theorem 3.1.]{Bor99}, this determines a precise collection of relations between the cycles $[\widetilde{\NL}_{g,d}]$. Conjecture \ref{conj: modularity} is similar in shape to several results about special cycles on orthogonal Shimura varieties in the Kudla-Millson program \cite{Bru02a, Bru02b, BLMM17}. To our knowledge, it is not clear that Conjecture \ref{conj: modularity} follows from these results.

\begin{rmk}
    In \cite{GLS24}, the authors show that the images of the cycles $\widetilde{\NL}_{g,d}$ under the map
    $$
    \widetilde{\NL}_{g,d} \longrightarrow\Ac_1 \times\Ac_{g-1}
    $$
    that remembers the domain and codomain of the homomorphism vanish in cohomology, but give non-zero modular forms after adding some level structure. This does not contradict our result, since $\Ac_1$ is not proper.
\end{rmk}

\subsection{The homomorphism property and Gromov-Witten theory of the universal elliptic curve}

If $\alpha$ is a tautological class, then for any $\beta$, 
\begin{equation}\label{eqn: semihomomorphism}
    \alpha \cdot\taut(\beta) = \taut(\alpha \cdot \beta)\, .
\end{equation}
Indeed, if $\Lambda$ is a tautological class,
\begin{align*}
    \int_{\overline{\Ac}_g} \taut(\alpha \cdot \beta) \cdot \Lambda \cdot \lambda_g &= \int_{\overline{\Ac}_g}  \beta\cdot \alpha \cdot \Lambda \cdot \lambda_g = \int_{\overline{\Ac}_g} \taut(\beta) \cdot \alpha \cdot \Lambda\cdot \lambda_g\, .
\end{align*}

\begin{defn}\label{defn: homomorphism property}
    Let $\alpha, \beta \in \CH^*(\Ac_g)$. We say that the pair $(\alpha, \beta)$ has the \emph{homomorphism property} if
    \begin{equation} \label{eqn:homomorphismproperty}
        \taut(\alpha \cdot \beta) = \taut(\alpha) \cdot \taut(\beta)\, .
    \end{equation}
\end{defn}

We do not know if the tautological projection is a ring homomorphism in general, but we find that this is the case in several situations
\begin{thm}\label{thm: homomorphism property}
    If $\alpha, \beta \in \CH^*(\Ac_g)$ are cycles supported on the Noether-Lefschetz locus given by non-simple abelian varieties, then both sides of equation \eqref{eqn:homomorphismproperty} vanish, so $(\alpha, \beta)$ has the homomorphism property.
\end{thm}

The authors in \cite{COP} define extensions of the tautological ring of van der Geer that include Noether-Lefschetz cycles. For instance, they consider the ring $\mathsf{R}^*_{\mathsf{pr}}(\Ac_g)$ generated by the pushforwards to $\Ac_g$ of classes in
$$
\Ac_{g_1} \times \ldots \times \Ac_{g_l}
$$
that can be written in terms of the $\lambda$ classes on each factor, or more generally the ring generated by all Noether-Lefschetz loci and the pushforwards of their tautological classes. Theorem \ref{thm: homomorphism property} shows that the tautological projection
$$
\taut : \mathsf{R}_{\mathsf{pr}}^*(\Ac_g) \longrightarrow \mathsf{R}^*(\Ac_g)
$$
is a ring homomorphism.

\begin{question}\label{conj: ring homomorphism}
    For which natural subrings of $\CH^*(\Ac_g)$ is the tautological projection a ring homomorphism?
\end{question}

Faber and Pandharipande introduced in \cite{FP03} the $\lambda_g$-pairing on the moduli space of curves of compact type:
\begin{equation}\label{eqn: Lg pairing curves}
    \CH^k(\Mc_g^{ct}) \times \CH^{2g-3-k}(\Mc_g^{ct}) \longrightarrow \CH^{2g-3}(\Mc_g^{ct}) \longrightarrow \Q\, ,
\end{equation}
given by $(\alpha, \beta) \mapsto \int_{\overline{\Mc}_g} \bar{\alpha}\cdot \bar{\beta}\lambda_g$, where $\lambda_g$ is the top Chern class of the Hodge bundle on $\overline{\Mc}_g$.

The Torelli morphism $\Tor : \Mc_g^{ct} \longrightarrow \Ac_g $ is proper, and its image is the locus of Jacobians, $\Jc_g$.

We want to investigate the homomorphism property for the pair $(\alpha , [\Jc_g])$, where $\alpha \in \CH^*(\Ac_g)$. By \cite{Nam76}, the Torelli morphism extends to a morphism between $\overline{\Mc}_g$ and some toroidal compactification of $\Ac_g$ while still preserving the Hodge bundle, and so $(\alpha, [\Jc_g])$ has the homomorphism property if and only if $\Tor^*(\alpha - \taut(\alpha))$ pairs to $0$ with any $\lambda$ class in $\Mc_g^{ct}$. A more refined result is given in Proposition \ref{prop: homomorphism and Gorenstein}.

Consider the universal family of elliptic curves with a zero section:
$$
\begin{tikzcd}
\overline{\Ec} \arrow[d, "\pi"']                  \\
{\overline{\Mc}_{1,1}} \arrow[u, "s"', bend right]
\end{tikzcd}.
$$
The stack $\overline{\Mc}_{g,1}(\pi,d)$ is the moduli space of log-stable maps of degree $d$ to the fibers of $\pi$ from a pointed curve of genus $g$. It is a proper Deligne--Mumford stack, is equipped with maps
$$
\mathsf{ev} : \overline{\Mc}_{g,1}(\pi,d) \longrightarrow \overline{\Ec}\, , \qquad \mathsf{ft} : \overline{\Mc}_{g,1}(\pi,d) \longrightarrow \overline{\Mc}_{g,1}\, ,
$$
and has virtual class of dimension $2g$. The  cotangent line at the marking of the domain curve defines a line bundle on $\overline{\Mc}_{g,1}(\pi,d)$, whose first Chern class is
$$
\psi \in \CH^1(\overline{\Mc}_{g,1}(\pi,d))\, .
$$
We consider Gromov-Witten invariants defined from this data:
$$
 \langle \tau_i(s)\Lambda \rangle_{g, d}^\pi := \int_{[\overline{\Mc}_{g,1}(\pi, d)]^{vir}}\psi^i \mathsf{ev}^*[s]\mathsf{ft}^*\Lambda\, ,
$$
where $\Lambda \in \CH^*(\overline{\Mc}_{g,1})$ and $i \geq 1$.

If $(C,p)$ is a marked curve of compact type, there is an Abel-Jacobi map $C \longrightarrow\Jac(C)$. If $\Jac(C)$ has an elliptic subgroup, the composition of the Abel-Jacobi map and the dual of the inclusion of $E$ gives rise to a morphism $(C,p) \longrightarrow(E,0)$ and every morphism from $C$ to an elliptic curve arises in this way. This suggests that $\mathrm{Tor}_1^{-1}(\NL_{g,d})$ is connected with the Gromov-Witten invariants arising from $\overline{\Mc}_{g,1}(\pi,d)$, where
$$
\Tor_1 : \Mc_{g,1}^{ct} \longrightarrow \Ac_g
$$
sends a pointed curve to its Jacobian. We will make this precise in Proposition \ref{prop: Tor to stable 1}, with an ingredient being the recent work of Greer and Lian (more concretely, \cite[Theorem 2]{GL24}) to identify $\mathrm{Tor}_1^*([\NL_{g,d}])$ with the virtual class of a component of $\mathsf{ft}^{-1}(\Mc_{g,1}^{ct})$.

If $([\NL_{g,d}], [\Jc_g])$ has the homomorphism property, we can replace
$$
[\overline{\Mc}_{g,1}(\pi,d)]^{vir} \cdot \mathsf{ev}^*[q]\lambda_g
$$ 
with $\Tor_{1}^*(\taut([\widetilde{\NL}_{g,d}])) \lambda_g$ for integral evaluations, thus arriving at the following conjecture:
\begin{conj}\label{conj: hilb C2}
    For all $g \geq 2$, $d \geq 1$, and $\Lambda \in \mathsf{R}^*(\overline{\Mc}_g)$ is an insertion divisible by $\lambda_g$ then 
    $$
    \langle \tau_i(s)\Lambda \rangle_{g, d}^\pi := \frac{g \sigma_{2g-1}(d)}{6|B_{2g}|}  
 \int_{\overline{\Mc}_{g,1}} \psi^i\lambda_{g-1} \Lambda\, ,
    $$
    where $B_{i}$ is the $i$-th Bernouilli number and $\sigma_{k}(n) = \sum_{m \mid n} m^k$.
\end{conj}
A particular case is
    \begin{equation}\label{eqn: Lg-2 evaluation}
    \langle \tau_1(s) \lambda_g \lambda_{g-2}\rangle_{g,d}^{\pi} = \frac{|B_{2g-2}|}{24(2g-2)!}\sigma_{2g-1}(d)\, ,
    \end{equation}
by \cite[Theorem 4]{FP00}. Notice that this formula has the same shape as the corresponding Gromov-Witten invariants of an elliptic curve, see \cite[Proposition 4.4.6]{P08}. The following result summarizes the connections:
\begin{thm}\label{thm: equivalent conj}
    Let $g \geq 2$. The following are equivalent:
    \begin{enumerate}
        \item[$a)$] Formula \eqref{eqn: Lg-2 evaluation} holds for all $d$.
        \item[$b)$] Conjecture \ref{conj: hilb C2} holds.
        \item[$c)$] For all $d$, the pair $([\NL_{g,d}], [\Jc_g])$ has the homomorphism property.
        \item [$d)$] For all $d$,
        $$
        \Tor^*([\NL_{g,d}] - \taut([\NL_{g,d}])) \in \CH^{g-1}(\Mc_g^{ct})
        $$
        pairs to $0$ with any tautological class in $\mathsf{R}^{g-2}(\Mc_g^{ct})$ under the $\lambda_g$ pairing \eqref{eqn: Lg pairing curves}.
    \end{enumerate}
\end{thm}

\begin{rmk}
    In a subsequent paper \cite{IPT24} with Pandharipande and Tseng, we will show that formula \eqref{eqn: Lg-2 evaluation} holds. This shows that all the statements in Theorem \ref{thm: equivalent conj} hold. We will also explain the connections to the quantum cohomology of the Hilbert scheme of points in $\C^2$. In particular, all the statements in Theorem \ref{thm: equivalent conj} are also equivalent to a conjecture by Pandharipande and Tseng on the simplest non-trivial genus 1 Gromov-Witten invariant of $\mathrm{Hilb}(\C^2)$, namely,
    $$
    \langle (2, 1^{d-2})\rangle_{g=1,\beta_n}^{\mathrm{Hilb}^d(\C^2)}\, .
    $$
\end{rmk}


\subsection{Plan of the paper}

In Section \ref{section: 2}, we introduce the morphisms $\phi$ and $\varphi$, which will be useful to deal with the $\NL_{g,\delta}$, and explain how they extend to the toroidal compactifications of the various moduli spaces. In Section \ref{section: 3}, we use these morphisms to prove Theorem \ref{thm: projection of NLD} up to the determination of the constant $c$, which is computed in Section \ref{section: 7}. We also show that $\taut(\alpha \cdot \beta)=0$ for any two cycles supported on the $\NL_{g,\delta}$, which is the major step in the proof of Theorem \ref{thm: homomorphism property}. In Section \ref{section: 4}, we discuss in more detail the Hecke operators and prove some of their properties. We also prove Theorem \ref{thm: not tautological}. In Section \ref{section: 5}, we finish the proof of Theorem \ref{thm: homomorphism property} and explain the connections of the homomorphism property to the kernel of the $\lambda_g$-pairing on the moduli space of curves. In Section \ref{section: 6}, we explain the connections to Gromov-Witten theory and the proof of Theorem \ref{thm: equivalent conj}.

Throughout the paper we work over the complex numbers, although all the results hold over any algebraically closed field of characteristic $0$.

\subsection{Acknowledgements}
I would like to thank Rahul Pandharipande for his patience, for proposing me this area of research and for all his support and guidance during the years of preparation of this article. I am very grateful to Gerard van der Geer for sharing his wisdom on the whole topic of moduli spaces of abelian varieties, and for suggesting to look at Hecke operators. I would also like to thank Younghan Bae, Samir Canning, Alessio Cela, Lycka Drakengren, Jeremy Feusi, Samouil Molcho, and Johannes Schmitt for the many conversations at ETH Zürich about abelian varieties, curves and intersection theory, and for comments on drafts of this paper, and the anonymous referee for the constructive comments and improvements. I am grateful to Qizheng Yin for his remarks on how to compute volumes of Hecke operators. Carl Lian and François Greer have been very helpful in explaining their work and conjectures.

The author was supported by ERC-2017-AdG-786580-MACI and SNF-200020-219369.

The project received funding from the European Research Council (ERC) under the European Union Horizon 2020 research and innovation programme (grant agreement 786580).

\newpage
\section{The morphisms \texorpdfstring{$\phi_\delta$}{phi d}, \texorpdfstring{$\pi_{\delta}$}{pi d} and \texorpdfstring{$\varphi_{g,\delta}$}{phi d}}\label{section: 2}

\subsection{Higher degree polarizations}\label{section: 21}
If $\theta : X \longrightarrow X^\vee$ is a polarization on an abelian variety, $\ker(\theta)$ carries a symplectic form, and given a Lagrangian subgroup $H$ of $\ker(\theta)$, there is a unique principal polarization $\theta_H$ on $X/H$.

If $X$ is an abelian variety with a polarization $\theta$ of type $\delta$, then a \emph{level structure} for $\theta$ is a symplectic isomorphism of groups
$$
    f: K(\delta) \stackrel{\cong}{\longrightarrow}\ker (\theta)\, ,
$$
where $K(\delta):=(\Z/d_1 \times \ldots \times \Z/d_g)^2$ carries a canonical symplectic form.
The moduli stack of triples $(X, \theta ,f)$ as above is denoted by $\Ac_{g,\delta}^\mathsf{lev}$ and it is a smooth Deligne-Mumford stack of dimension $\binom{g+1}{2}$.
\begin{defn}
    The morphism
    $$
    \phi_\delta : \Ac_{g,\delta}^\mathsf{lev} \longrightarrow \Ac_g
    $$
    is obtained by sending $(X,\theta ,f)$ to $(X/f(\{0\} \times\Z^g/\delta ), \theta_{f(\Z^g/\delta \times \{0\})})$.\\
    There is also a morphism
    $$
    \pi_\delta : \Ac_{g,\delta}^\mathsf{lev}\longrightarrow\Ac_{g,\delta}
    $$
    that forgets the level structure.
\end{defn}

Both $\phi_\delta$ and $\pi_\delta$ are finite morphisms, and $\phi_{\delta}^*\mathbb E_g = \pi_\delta^*\mathbb E_{g,\delta}$, since the Hodge bundle is preserved under isogenies.

\begin{lem}\label{lem: alpha over the complex}
    Consider the following groups:
    \begin{align*}
    \Grm_{\delta} &= \left\lbrace M \in \Sp_{2g}(\Q)\left|\, M^t \begin{pmatrix}
        \mathbbm{1}_g&0\\0&\delta
    \end{pmatrix}\Z^{2g} \subseteq \begin{pmatrix}
        \mathbbm{1}_g&0\\0&\delta
    \end{pmatrix}\Z^{2g}\right.\right \rbrace\, , \\
    \Grm_\delta [\delta] &= \left\lbrace M \in \Sp_{2g}( \Z) \mid M = \mathbbm{1}_{2g} + \begin{pmatrix}
        \delta&0 \\
        0 &\mathbbm{1}_g
    \end{pmatrix}N \begin{pmatrix}
        \mathbbm{1}_g & 0\\
        0 &\delta
    \end{pmatrix}\text{ for some }N \in \mathrm{M}_{2g}(\Z) \right\rbrace\, .
    \end{align*}
    Over the complex numbers, one can realize $\Ac_{g,\delta}$ and $\Ac_{g,\delta}^\mathsf{lev}$ as the quotients $[\Hc_g / \Grm_{\delta}]$ and $[\Hc_g /\Grm_\delta[\delta]]$ respectively, so that, 
    under this identification, $\pi_\delta$ and $\phi_\delta$ are induced by the inclusions $\Grm_{\delta}[\delta] \subset \Grm_{\delta}$ and $\Grm_{\delta}[\delta] \subset \mathrm{Sp}(2g,\Z)$.
\end{lem}
\begin{proof}
    The stack $\Ac_{g,\delta}^\mathsf{lev}$ is explained in \cite[Section 8.3.1.]{BL04}, and from the comments in Section 8.2., it is clear that, in their language,
    $\Ac_{g,\delta}^\mathsf{lev} = [\Hc_g / \sigma_{\delta}(\Gamma_\delta(\delta))]$, where
    \begin{equation*}
        \begin{alignedat}{2}
        \sigma_{\delta} : \left\lbrace M \in \mathrm{M}_{2g}(\Z) \left| \, M \begin{pmatrix}
        0&\delta\\
        -\delta &0
    \end{pmatrix}M^t = \begin{pmatrix}
        0&\delta\\
        -\delta &0
    \end{pmatrix}\right.\right\rbrace &&\longrightarrow&\Sp_{2g}(\Q)\\
    M &&\mapsto &\begin{pmatrix}
        \mathbbm{1}_g&0\\
        0&\delta
    \end{pmatrix} M \begin{pmatrix}
        \mathbbm{1}_g&0\\
        0&\delta
    \end{pmatrix}^{-1}\, .
    \end{alignedat}
    \end{equation*}
    One can easily check that $\sigma_{\delta}(\Gamma_\delta(\delta))$ is precisely the group $\mathrm{G}_{g}[\delta]$. Under the aforementioned identification,
    $\tau \in \Hc_g$ corresponds to the triplet
    $$
    \left(\frac{\C^g}{( \tau \, |\, \delta)}, \mathrm{Im}(\tau)^{-1}, (\delta^{-1}\cdot \tau \, | \, \mathbbm{1}_g)\right) \in \Ac_{g, \delta}^{\mathsf{lev}}\, .
    $$
    and with this description, $\phi_\delta$ quotients out the second factor of the kernel of the polarization, so $\phi_\delta$ sends the polarized abelian variety $(\C^g /(\tau\, |\, \delta), \mathrm{Im}(\tau)^{-1})$ to the principally polarized abelian variety $(\C^g /(\tau \, | \,\mathbbm{1}_g), \mathrm{Im}(\tau)^{-1})$. Therefore, $\phi_\delta$ is induced by the identity on $\Hc_g$. For the map $\pi_\delta$ this follows from the construction in \cite[Section 8.3.1.]{BL04}.
    \end{proof}

\begin{proof}[Proof of Theorem \ref{thm: tautological rings}]
Theorem \ref{thm: tautological rings} holds when $\delta = (1, \ldots , 1)$ by \cite[Theorem 3]{vdG99}. For an arbitrary $\delta$, consider the diagram
    $$
    \begin{tikzcd}
                 & {\Ac_{g,\delta}^\mathsf{lev}} \arrow[ld, "\pi_\delta"'] \arrow[rd, "\phi_\delta"] &       \\
{\Ac_{g,\delta}} &                                                                              & \Ac_g
\end{tikzcd}\, ,
    $$
    which gives rise to an isomorphism
    $$
    \left(\frac{1}{\deg(\pi_\delta)}\pi_{\delta, *}\right) \circ \phi_\delta^* : \mathsf{R}^*(\Ac_g) \longrightarrow \mathsf{R}^*(\Ac_{g,\delta})
    $$
    sending $c_i(\mathbb E_{g})$ to $c_i(\mathbb E_{g,\delta})$, because the Hodge bundle is preserved by isogenies and so $\pi_\delta^* \mathbb E_{g,\delta} = \phi_\delta^*\mathbb E_g$.
\end{proof}

\subsection{Non-simple abelian varieties}\label{section: 22}

In this section we introduce the tools to work with the locus of principally polarized abelian varieties having an abelian subvariety. We will see that this is a countable union of substacks of $\Ac_g$, indexed by the type of the induced polarization on the subvariety.

If $\delta =(d_1, \ldots , d_u)$ is a polarization type with $u\leq g/2$, its \textit{complementary type} is
$$
\widetilde{\delta} = (\underbrace{1, \ldots ,1}_{g-2u \text{ times}},d_1, \ldots , d_u)\, ,
$$
and the \textit{double type} is
$$
\delta^2 = (\underbrace{1, \ldots ,1}_{g-2u \text{ times}},d_1,d_1, \ldots ,d_u, d_u)\, .
$$
If $(X, \theta)$ is a principally polarized abelian variety of dimension $g$ and $Y$ is a $u$-dimensional abelian subvariety, then
$$
Z:= \ker \left(Y \stackrel{\theta}{\longrightarrow} X^\vee \longrightarrow Y^\vee\right)^{\circ}
$$
is the \textit{complementary subvariety} of $Y$ It has dimension $g-u$ and by \cite[Lemma 12.1]{BL04}, if $\theta|_{ Y}$ is of type $(d_1, \ldots , d_u)$ then $\theta|_{ Z}$ is of type $\widetilde{\delta}$.

Conversely, given polarized abelian varieties
$$
(X, \theta_X) \in \Ac_{u,\delta}\hspace{6pt}\text{ and } \hspace{6pt}(Y, \theta_Y) \in \Ac_{g-u, \widetilde{\delta}}\, ,
$$
we can construct principally polarized abelian varieties $X$ containing $Y$ and $Z$ as complementary subvarieties with the prescribed polarization. This construction is due to Debarre \cite{D88}, and is also explained with detail in \cite[Section 3]{A16}. We summarize the construction in the following lemma, which works in families.
\begin{lem}\label{lem: graphs and polarisations}
    With notation as above, for any antisymplectic\footnote{If $(A, \omega_A)$ and $(B, \omega_B)$ are two abelian groups with an alternating bilinear form, a linear map $f:A \longrightarrow B$ is \emph{antisymplectic} if $\omega_B(f(a_1), f(a_2)) = - \omega_A(a_1, a_2)$ for all $a_1, a_2 \in A$.} isomorphism
    $$
    p : \ker (\theta_Y) \longrightarrow\ker(\theta_Z)\, ,
    $$
    the abelian variety
    $$
    \frac{Y \times Z}{\mathrm{graph}(p)}
    $$
    has a canonical principal polarization $\theta_p$ and contains $(Y, \theta_Y)$, $(Z,\theta_Z)$ as complementary subvarieties. Moreover, the isomorphism type of
    $$
    \left(\frac{Y \times Z}{\mathrm{graph}(p)}, \theta_p\right)
    $$
    does not depend on $p$, and all principally polarized abelian varieties having $Y$ and $Z$ as complementary subvarieties arise this way.
\end{lem}

\begin{defn}
    If $u \leq g/2$ and $\delta = (d_1, \ldots ,d_u)$ is a polarization type, define the morphism
    $$
    \varphi_{g,\delta} :  \Ac_{u, \delta}^{\mathsf{lev}} \times \Ac_{g-u, \widetilde{\delta}}^\mathsf{lev} \longrightarrow \Ac_g
    $$
    sending the pair $((Y, \theta_Y , f_Y), (Z, \theta_Z , f_Z))$ to 
    $$
    \left(\frac{Y \times Z}{\mathrm{graph} (f_Z \circ r \circ f_Y^{-1})},\theta_{f_Z \circ r \circ f_Y^{-1}}\right),
    $$
     where
    $$
    r: K(\delta) \longrightarrow K(\widetilde{\delta}) = K(\delta)
    $$
     is the antisymplectic isomorphism that exchanges the factors $(\Z^g/\delta) \times 0$ and $0 \times (\Z^g/\delta)$.
     By Lemma \ref{lem: graphs and polarisations}, $\varphi_{g,\delta}$ is surjective onto the locus of principally polarized abelian varieties having a subvariety whose induced type is of type $\delta$. In particular, this locus is irreducible.
\end{defn}

A symplectic automorphism $ s: K(\delta) \longrightarrow K(\delta)$ acts on $\Ac_{u, \delta}^\mathsf{lev} \times \Ac_{g-u, \tilde{\delta}}^\mathsf{lev}$ via
$$
s.((Y, \theta_Y , f_Y), (Z, \theta_Z , f_Z)) = ((Y, \theta_Y , f_Y\circ s^{-1}), (Z, \theta_Z , f_Z\circ s^{-1}))\, ,
$$
and $\varphi_{g, \delta}$ is invariant under this action. Let $\Ac_{u, g-u, \delta}$ be the quotient by the action of $\Sp(K(\delta))$ just described. It parametrizes triplets $(X,Y,Z, \theta)$ where $(X, \theta)$ is a principally polarized abelian variety of dimension $g$, and $Y$, $Z$ are complementary subvarieties of a principally polarized abelian variety $X$, of dimensions $u$ and $g-u$, and the type of the induced polarizations $\theta |_{Y}$, $\theta |_{Z}$ is $\delta$ and $\widetilde{\delta}$, respectively.

There is a commutative diagram
$$
\begin{tikzcd}
	{\Ac_{u, \delta}^\mathsf{lev} \times \Ac_{g-u, \widetilde{\delta}}^\mathsf{lev}} \\
	& {\Ac_{u, g-u, \delta}} & {\Ac_{g}} \\
	{\Ac_{u, \delta} \times \Ac_{g-u, \widetilde{\delta}}}
	\arrow["\tau", from=1-1, to=2-2]
	\arrow["{\varphi_{g,\delta}}", bend left=20, from=1-1, to=2-3]
	\arrow["{\pi_{\delta} \times \pi_{\delta}}"', from=1-1, to=3-1]
	\arrow["{\varphi'}", from=2-2, to=2-3]
	\arrow[from=2-2, to=3-1]
\end{tikzcd} \, ,
$$
and $\tau$ is a principal $\Sp(K(\delta))$-bundle, so $\Ac_{u, g-u, \delta}$ is a smooth Deligne-Mumford stack.

Let $\mathbb E_{u,\delta}$ and $\mathbb E_{g-u, \widetilde{\delta}}$ be the pullbacks of the Hodge bundles from $\Ac_{u, \delta}$ and $\Ac_{g-u, \widetilde{\delta}}$. Then the tangent bundle of $\Ac_{u, g-u, \delta}$ is $\operatorname{Sym}^2\mathbb E_{u, \delta}^\vee \boxplus \operatorname{Sym}^2\mathbb E_{g-u, \tilde{\delta}}^\vee$.

Since the morphism $\varphi'$ is unramified between smooth DM stacks, it is, \'etale locally, a regular immersion by \cite[Corollary 18.4.7.]{EGA4}. The pullback of the Hodge bundle on $\Ac_g$ via $\varphi'$ splits as $\mathbb E_{u, \delta} \boxtimes \mathbb E_{g-u, \tilde{\delta}}$, so the normal bundle to $\varphi'$ equals
$$
N_{\varphi'} = \operatorname{Sym}^2(\mathbb E_{u, \delta}^\vee \boxtimes \mathbb E_{g-u, \tilde{\delta}}^\vee) - \operatorname{Sym}^2\mathbb E_{u, \delta}^\vee \boxplus \operatorname{Sym}^2\mathbb E_{g-u, \tilde{\delta}}^\vee =  \mathbb E_{u, \delta}^\vee \boxtimes \mathbb E_{g-u, \tilde{\delta}}^\vee\, .
$$
\begin{defn}
    $\NL_{g,\delta}$ is the stack theoretic image of  $\varphi'$.
\end{defn}
In particular, since generically abelian varieties are simple, $\varphi'$ has degree $1$, so
$$
[\NL_{g,\delta}] = \varphi'_*(\mathbf{1}) = \frac{1}{|\Sp(K(\delta))|} \varphi_{g,\delta, *}(\mathbf{1})\in \CH^{u(g-u)}(\Ac_g)\, .
$$

\begin{rmk}
    When $u = g/2$, $\varphi'$ has degree $2$ because one can exchange the role of $Y$ and $Z$ so $[\NL_{g,\delta}]$ should be $\frac{1}{2}\varphi'_*(\mathbf{1})$ in principle. However, it is more convenient to consider
    $$
    [\NL_{g,\delta}] := \varphi'_*(\mathbf{1})\, .
    $$
    This is consistent with \cite{COP, CMOP24}, where the cycle classes $[\Ac_{g/2} \times \Ac_{g/2}]$ are the pushforward of $1$ along $\Ac_{g/2} \times \Ac_{g/2} \longrightarrow\Ac_g$.
\end{rmk}

\subsection{Extensions to the boundary}\label{section: 23}

Mumford and his collaborators introduce in \cite{AMRT10} various toroidal compactifications of the quotient of a hermitian symmetric domain $\Dc$ (such as the Siegel upper half space) by an arithmetic subgroup $\Gamma$ of its group of automorphisms. These compactifications depend on a choice of a $\Gamma$-admissible subdivision $\Sigma$ of a certain cone. The resulting proper space is denoted by
$$
\overline{[\Dc/\Gamma]}^\Sigma\, .
$$
If $\Dc = G/K$, where $K$ is a maximal compact  subgroup, complex representations of $K$ give rise to \emph{automorphic bundles} on $[\Dc/\Gamma]$, such as the Hodge bundle (which arises as the defining representation of $U(g) \subseteq \Sp(2g, \mathbb R)$ on $\mathbb C^g$). These bundles have a canonical extension to any toroidal compactification \cite{M77}. We denote the extension of the Hodge bundle by $\mathbb E_g$.

The morphisms $\pi_\delta$, $\phi_\delta$ and $\varphi_{g, \delta}$ extend to arbitrarily fine toroidal compactifications of the moduli spaces. This follows from the results in \cite{Ha89}, so it might be clear to experts, but we explain it in detail.

Given two arithmetic subgroups $\Gamma < \Gamma'$ for a hermitian symmetric domain $\Dc$, any subdivision $\Sigma$ which is $\Gamma'$ admissible is also $\Gamma$-admissible and the map
$$
\overline{[\Dc/\Gamma]}^{\Sigma} \longrightarrow \overline{[\Dc/\Gamma']}^{\Sigma}
$$
is well-defined. The morphisms $\pi_\delta$ and $\phi_\delta$ are quotients by a finite subgroup, so their extendability is proven.

For $\varphi_{g,\delta}$, note that we can factor it as
$$
\Ac_{u, \delta}^\mathsf{lev} \times \Ac_{g-u, \widetilde{\delta}}^\mathsf{lev}\stackrel{\iota}{\longrightarrow} \Ac_{g, \delta^2}^\mathsf{lev} \stackrel{\kappa}{\longrightarrow} \Ac_{g}
$$
where $\iota$ takes two polarized abelian varieties to their product, which is polarized by the product of the two polarizations, and sends two bases of the kernels of the polarizations to their union, and $\kappa$ sends an abelian variety to its quotient by the antidiagonal of the identity $K(\delta) \longrightarrow K(\delta)$, seen as a subset of $K(\delta^2)$. Then, $\kappa$ is just a quotient by a congruence subgroup as before,  and with the identifications as in Lemma \ref{lem: alpha over the complex}, $\iota$ is induced by the map
\begin{equation*}
\begin{alignedat}{2}
    \Hc_u \times \Hc_{g-u} &\longrightarrow && \quad \,\,\Hc_g\\
    (\tau', \tau'')\quad &\mapsto &&\begin{pmatrix}
        \tau' &0\\
        0&\tau''
    \end{pmatrix}
\end{alignedat}\, ,
\end{equation*}
which is a morphism of hermitian symmetric domains, and so $\iota$ extends to arbitrary toroidal compactifications by the explanation in \cite[Section 3.3.]{Ha89}. The operations of pullback an automorphic bundle by (an extension of) a morphism of hermitian symmetric domain of an automorphic bundle and taking canonical extensions commute by \cite[Theorem 4.2.]{Ha89} and so,
$$
\overline{\varphi}_{g, \delta}^* \mathbb E_g = \mathbb E_u \boxplus \mathbb E_{g-u} \quad \overline{\pi}_{\delta}^* \mathbb E_g = \mathbb E_g \quad \overline{\phi}^*_\delta \mathbb E_g = \mathbb E_g \, .
$$

\section{Intersection properties of \texorpdfstring{$\NL_{g,\delta}$}{NL gd}}\label{section: 3}
\subsection{Tautological projection}\label{section: 31}
Recall that the tautological projection $\taut(\alpha)$ is the unique tautological class such that
    $$
    \int_{\overline{\Ac}_g} \overline{\alpha} \cdot \Lambda \cdot \lambda_g = \int_{\overline{\Ac}_g} \taut(\alpha) \cdot \Lambda \cdot \lambda_g
    $$
    for any tautological class $\Lambda \in \mathsf{R}^*(\Ac_g)$, where $\overline{\alpha}$ is an extension of $\alpha$ to $\overline{\Ac}_g$.

\begin{proof}[Proof of Theorem \ref{thm: projection of NLD}]

Consider extensions of the morphisms $\phi_\delta \times \phi_{\tilde{\delta}}$ and $\varphi_{g, \delta}$ as discussed in in Section \ref{section: 23}, so that there is a commutative diagram
$$
\begin{tikzcd}
	{\overline{\Ac_{u,\delta}^\mathsf{lev} \times \Ac_{g-u, \tilde{\delta}}^\mathsf{lev}}} && {\overline{\Ac}_g}^{\Sigma_1} \\
	& {\Ac_{u,\delta}^\mathsf{lev} \times \Ac_{g-u, \widetilde{\delta}}^\mathsf{lev}} & {\Ac_{g}} \\
	& {\Ac_u \times \Ac_{g-u}} & {\Ac_g} \\
	{\overline{\Ac_u \times \Ac_{g-u}}} && {\overline{\Ac}_g}^{\Sigma_2}
	\arrow["{\overline{\varphi}_{g,\delta}}", from=1-1, to=1-3]
	\arrow["{\overline{\phi_{\delta}\times \phi_{ \tilde{\delta}}}}"', from=1-1, to=4-1]
	\arrow[hook', from=2-2, to=1-1]
	\arrow["{\varphi_{g,\delta}}", from=2-2, to=2-3]
	\arrow["{\phi_\delta \times \phi_{\tilde{\delta}}}", from=2-2, to=3-2]
	\arrow[hook', from=2-3, to=1-3]
	\arrow["{\varphi_{g,(1^u)}}"', from=3-2, to=3-3]
	\arrow[hook', from=3-2, to=4-1]
	\arrow[hook', from=3-3, to=4-3]
	\arrow["{\overline{\varphi}_{g,(1^u)}}"', from=4-1, to=4-3]
\end{tikzcd} \, ,
$$
and the pullbacks of the Hodge bundle to $\overline{\Ac_{u,\delta}^\mathsf{lev} \times \Ac_{g-u, \tilde{\delta}}^\mathsf{lev}}$ agree. Then $\overline{\varphi}_{g, \delta, *}(1)$ is an extension of $[\NL_{g,\delta}]$. If $\Lambda$ is a tautological class, and we write $\Lambda \cdot \lambda_{g} = P(c_1(\mathbb E_g), \ldots , c_g(\mathbb E_g))$,
    \begin{align*}
        \int_{\overline{\Ac}^{\Sigma_1}_g} \overline{[\NL_{g, \delta}]} \cdot \Lambda \cdot \lambda_g & = \frac{1}{\deg (\overline{\varphi}_{g, \delta})}\int_{\overline{\Ac}_g} \overline{\varphi}_{g, \delta, *}(1) \cdot P(c_i(\mathbb E_g))\\
        &= \frac{1}{\deg (\overline{\varphi}_{g, \delta})}\int_{\overline{\Ac_{u,\delta}^\mathsf{lev} \times \Ac_{g-u, \tilde{\delta}}^\mathsf{lev}}} P(c_i(\overline{\varphi}_{g, \delta}^*\mathbb E_g))\\
        &=  \frac{1}{\deg (\overline{\varphi}_{g, \delta})}\int_{\overline{\Ac_{u,\delta}^\mathsf{lev} \times \Ac_{g-u, \tilde{\delta}}^\mathsf{lev}}} P(c_i(\overline{\phi_{\delta}\times \phi_{\tilde{\delta}}}^*(\mathbb E_u \boxplus \mathbb E_{g-u})))\\
        &=\frac{\deg (\overline{\phi_{\delta}\times \phi_{ \tilde{\delta}})}}{\deg (\overline{\varphi}_{g, \delta})}\int_{\overline{\Ac_u \times \Ac_{g-u}}} P(c_i(\mathbb E_u \boxplus \mathbb E_{g-u}))\\
        &=\frac{\deg (\overline{\phi_{\delta}\times \phi_{ \tilde{\delta}})}}{\deg (\overline{\varphi}_{g, \delta})} \int_{\overline{\Ac_u \times \Ac_{g-u}}} P(c_i(\overline{\varphi}_{g, (1^u)}^*\mathbb E_g))\\
        & = \frac{\deg (\overline{\phi_{\delta}\times \phi_{ \tilde{\delta}})}}{\deg (\overline{\varphi}_{g, \delta})} \int_{\overline{\Ac}^{\Sigma_2}_g} \overline{[\Ac_u \times \Ac_{g-u}]} \cdot \Lambda \cdot \lambda_g\\
        &= \frac{\deg (\overline{\phi_{\delta}\times \phi_{ \tilde{\delta}})}}{\deg (\overline{\varphi}_{g, \delta})} \int_{\overline{\Ac}^{\Sigma_2}_g} \taut([\Ac_u \times \Ac_{g-u}]) \cdot \Lambda \cdot \lambda_g \\
        &= \int_{\overline{\Ac}^{\Sigma_1}_g} \frac{\deg (\overline{\phi_{\delta}\times \phi_{ \tilde{\delta}})}}{\deg (\overline{\varphi}_{g, \delta})} \taut([\Ac_u \times \Ac_{g-u}]) \cdot \Lambda \cdot \lambda_g\, ,
    \end{align*}
    by several applications of the projection formula. Note also that
    $$
    \frac{\deg (\overline{\phi_{\delta}\times \phi_{ \tilde{\delta}})}}{\deg (\overline{\varphi}_{g, \delta})}  = \frac{\deg (\phi_{\delta})\cdot\deg( \phi_{\tilde{\delta}})}{\deg (\varphi_{g, \delta})} = \frac{\deg (\phi_{\delta})\cdot\deg( \phi_{\tilde{\delta}})}{|\Sp(K(\delta))|}\, .
    $$
    by Lemma \ref{lem: graphs and polarisations}, and these numbers are calculated in Propositions \ref{prop: degree aD} and \ref{prop: degree piD}.
\end{proof}

\begin{rmk}
    When $u=1$, one can prove that the tautological projections of all the $[\NL_{g,d}]$ are proportional to $\lambda_{g-1}$ in a different way:
    \begin{align*}
    \lambda_{g-1}\cdot [\NL_{g,d}] &= \frac{1}{\deg(\varphi_{g,d})}\varphi_{g,d,*}\varphi_{g,d}^*\lambda_{g-1} = \frac{1}{\deg(\varphi_{d,g})}\varphi_{g,d,*} c_{g-1}(\mathbb E_1 \boxplus \mathbb E_{g-1})\\
    &= \frac{1}{\deg(\varphi_{g,d})}\varphi_{g,d,*}(c_1(\mathbb E_1)c_{g-2}(\mathbb E_{g-1}) + c_{g-1}(\mathbb E_{g-1}))=0\, ,
    \end{align*}
    where we are using Theorem \ref{thm: tautological rings}. Therefore,
    \begin{align*}
    \lambda_{g-1} \taut(\NL_{g,d}) = \taut(\lambda_{g-1} \cdot [\NL_{g,d}])=0\, ,
    \end{align*}
    by formula \eqref{eqn: semihomomorphism}; but $\Q\lambda_{g-1}\subset \mathsf{R}^{g-1}(\Ac_{g})$ is precisely the kernel the map that sends $\alpha$ to $\alpha \cdot \lambda_{g-1}$. This was pointed out by R. Pandharipande.
\end{rmk}
\begin{proof}[Proof of Corollary \ref{cor: modular form}]
    Given a homomorphism $f: E \longrightarrow(A,\theta)$, let $d = \deg(f^*\theta)$ and $\hat{d} = d/|\ker (f)|$. This defines an elliptic subgroup $E'$ of $A$ such that $\theta_{E'}$ is of type $\hat{d}$. Conversely, given an elliptic subgroup $E'\leq A$, there are $\sigma_1(d/\hat{d})$ isogenies $E \longrightarrow E'$ of degree $d/\hat{d}$, where $\sigma_k(n) = \sum_{m \mid n} m^k$. Therefore,
    \begin{equation}\label{eqn: tilde to not tilde}
    [\widetilde{\NL}_{g,d}] = \sum_{\hat{d}\mid d}\sigma_1(d/\hat{d})[\NL_{g,\hat{d}}]\, ,
    \end{equation}
    and, by Theorem \ref{thm: projection of NLD},
    \begin{align*}
    \taut([\widetilde{\NL}_{g,d}]) &=\frac{g}{6|B_{2g}|}\sum_{\hat{d}\mid d}\sigma_{1}(d/\hat{d})d^{2g-1}\prod_{p \mid \hat{d}}(1-p^{2g-2}) \lambda_{g-1}\\
    &= \frac{g}{6|B_{2g}|}d\sum_{\hat{d}\mid d}\sigma_{-1}(d/\hat{d}) \hat{d}^{2g-2}\prod_{p \mid \hat{d}}(1-p^{2g-2})\lambda_{g-1}\, .
    \end{align*}
    The function $n \mapsto n^k \prod_{p \mid n} (1-p^{-k})$ is Jacobi's totient function $J_k$, and can also be written as a convolution $\operatorname{Id}_k * \mu$, where $\operatorname{Id}_k(n)=n^k$ and $\mu$ is the Möbius inversion function. By basic convolution identities, we have that
    $$
    d\cdot (\sigma_{-1}*J_{2g-2})(d)=d\cdot (\operatorname{Id}_{-1} * 1*\mu*\operatorname{Id}_{2g-2})(d) = d\cdot (\operatorname{Id}_{-1}*\operatorname{Id}_{2g-2})(d) = \sigma_{2g-1}(d)\, .
    $$
\end{proof}

\subsection{Intersections between \texorpdfstring{$\NL$}{NL} cycles}\label{section: 32}

In \cite{COP} the authors use the excess intersetion formula to show that $[\Ac_{g_1}\times \ldots \times \Ac_{g_l}] \cdot [\Ac_{h_1}\times \ldots \Ac_{h_k}]=0$. More generally, we study the intersection of any pair of cycles supported on the $\NL$ locus.
\begin{prop}\label{prop: taut of a product}
    Let $u\leq v \leq g/2$ and $\delta = (d_1, \ldots , d_u)$, $\varepsilon = (e_1, \ldots , e_v)$ be polarization types. Suppose that $\gamma, \eta$ are two cycles supported on $\NL_{g, \delta}$ and $\NL_{g, \varepsilon}$, respectively. Then
    $$
    \taut(\gamma \cdot \eta) =0\, .
    $$
\end{prop}
\begin{proof}
We form the Cartesian diagram
\begin{equation}\label{eqn: fiber product}
        \begin{tikzcd}
	\mathcal{Z} & \Ac_{v,g-v,\varepsilon} \\
	\Ac_{u, g-u,\delta} & {\Ac_g}
	\arrow[from=1-2, to=2-2]
	\arrow["\varphi_{g,\delta}",from=2-1, to=2-2]
	\arrow[from=1-1, to=1-2]
	\arrow[from=1-1, to=2-1]
	\arrow["\lrcorner"{anchor=center, pos=0.125}, draw=none, from=1-1, to=2-2]
\end{tikzcd}.
\end{equation}

$\mathcal Z$ parametrizes the data of an abelian variety with two pairs of complementary subgroups. The different incidences between these subgroups determine disconnected components of $\mathcal Z$. More concretely, let $X \sim X_1\times X_2 \sim Y_1 \times Y_2$, where $\dim (X_1)=u$, $\dim (Y_1) =v$ and $\sim$ denotes an isogeny. Up to isogeny, we have maps
$$
\pi_1, \pi_2 : X_1 \longrightarrow Y_1, Y_2
$$
and $\ker (\pi_1) \cap \ker (\pi_2)$ is a finite group, so
\begin{align*}
    X_1 &\sim \ker (\pi_1) \times \ker (\pi_2) \times E\, ,\\
    Y_1 &\sim \ker (\pi_2) \times E \times D_1\, ,\\
    Y_2 & \sim \ker(\pi_1) \times E \times D_2\, ,\\
    X_2 &\sim E \times D_1 \times D_2\, ,
\end{align*}
for some abelian varieties $E, D_1, D_2$. Note that the isogeny class of $E$ appears twice in the isogeny class of $X$. From the way we have constructed these varieties, it is clear that the dimensions of $E, D_1, D_2, \ker(\pi_1), \ker(\pi_2)$ stay constant if $X$ moves along a connected component of $\mathcal Z$. There are discrete invariants involving the types of the different polarizations induced on each of the varieties. The construction of $\varphi_{g,\delta}$ in Section \ref{section: 22} can be generalized, without having as much control on the degrees of the morphisms. After adding a large enough level structure, the abelian variety $X$, with its subvarieties $X_i, Y_i$, can be recovered from some polarized abelian varieties $X_i', Y_i'$ which are isogenous to $X_i$ and $Y_i$ respectively, with bounded isogeny degree, and with a level structure on $X_i'$ and $Y_i'$. 

Therefore, each connected component of $\mathcal Z$ has a finite morphism from
\begin{equation}\label{eqn: W definition}
\mathcal W := \Ac_{n_1,\mathsf{lev}}\times \Ac_{n_2,\mathsf{lev}}\times\Ac_{m_1,\mathsf{lev}}\times\Ac_{m_2,\mathsf{lev}}\times\Ac_{p,\mathsf{lev}}, ,
\end{equation}
where $n_i = \dim (\ker(\pi_i))$, $m_i = \dim(D_i)$ and $p = \dim(E)$, and $\Ac_{h,\mathsf{lev}}$ denotes a (connected) moduli space of polarized abelian varieties of dimension $h$ with some level structure\footnote{Here, \emph{level} has a more general meaning than the symplectic level structure. It could include the data of a basis of the $n$-torsion points of the variety, for some $n$. In any case, there is always a common roof between any two moduli spaces of abelian varieties of the same dimension with level structure.}, so that the maps to $\Ac_{v, g-v, \varepsilon}$  and $\Ac_{u, g-u, \delta}$ (which we will think about for the rest of the section as $\Ac_{v,\mathsf{lev}} \times \Ac_{g-v, \mathsf{lev}}$ and $\Ac_{u,\mathsf{lev}} \times \Ac_{g-u, \mathsf{lev}}$ respectively) factor as the diagonal on $\Ac_{p,\mathsf{lev}}$:
$$
\Delta : \mathcal W \longrightarrow \mathcal V:=\Ac_{n_1,\mathsf{lev}}\times \Ac_{n_2,\mathsf{lev}}\times\Ac_{m_1,\mathsf{lev}}\times\Ac_{m_2,\mathsf{lev}}\times\Ac_{p,\mathsf{lev}} \times \Ac_{p,\mathsf{lev}} 
$$
followed by multiplication maps
$$
\begin{alignedat}{2}
      \mathcal V & \stackrel{\mathrm{mult}_1}{\longrightarrow}  && \Ac_{v, \mathsf{lev}}\times \Ac_{g-v,\mathsf{lev}} \\
              (K_1, K_2, D_1, D_2, E_1, E_2) & \mapsto{} && (K_2\times D_1 \times E_1, K_1 \times D_2 \times E_2)
\end{alignedat}\, ,
$$

$$
\begin{alignedat}{2}
      \mathcal V & \stackrel{\mathrm{mult}_2}{\longrightarrow}  && \Ac_{u, \mathsf{lev}}\times \Ac_{g-u,\mathsf{lev}} \\
              (K_1, K_2, D_1, D_2, E_1, E_2) & \mapsto{} && (K_1\times K_2 \times E_1, D_1 \times D_2 \times E_2)
\end{alignedat}\, ,
$$
and possibly quotient maps such as $\phi_\delta$ or $\pi_\delta$ that forget level structure or send an abelian variety to the quotient of it by a finite group. We will omit these since they do not play a role in the computation.

We do this in detail when $\varepsilon = (1, \ldots , 1)$ and $u=1$, which is the only statement that we need for Theorem \ref{thm: equivalent conj}.
In this case, $\Zc$ parametrizes principally polarized abelian varieties $X = Y_1 \times Y_2$ that have an elliptic subgroup $X_1$ such that $\theta_X \cdot X_1 = d$. If $\pi_1(X_1)=0$ then $Y_2 \in \NL_{v,d}$ and the map
$$
\NL_{v,d} \times \Ac_{g-v} \longrightarrow\Zc
$$
is an isomorphism onto a connected component, and so we have a finite map
$$
\Ac_{v-1,(1,\ldots , 1,d)}^\mathsf{lev} \times \Ac_{1,(d)}^\mathsf{lev} \times \Ac_{g-v} \longrightarrow \Zc \, .
$$
Similarly if $\pi_2(X_1)=0$. Otherwise, $\pi_1$ and $\pi_2$ are non-zero maps. Define the discrete data $G_i = \ker(\pi_i)$, $d_i = |G_i|$ and $f_i = \pi_i(X_1)\cdot \theta_{X_i}$, so that
$$
G_1 \cap G_2 =\{0\}, \pi_i(X_1) = X_1/G_i \text{ and } d =d_1f_1 + d_2f_2 \, .
$$
Consider the modular curve
$$
\Ac_{1,G_1,G_2, f_1, f_2}^{\mathsf{lev}}
$$
that parametrizes elliptic curves $X_1'$ with a pair of subgroups $G_1', G_2'$, isomorphic to $G_1$ and $G_2$ as abstract groups that intersect trivially and a basis of the $f_1$-torsion points of $X_1'/G_1'$, and a basis of the $f_2$-torsion points of $X_1'/G_2'$. This modular curve is a finite cover of $\Ac_1$, and has a "digaonal" map
$$
\begin{alignedat}{2}
    \Delta : \Ac_{1,G_1,G_2, f_1, f_2}^{\mathsf{lev}} & \longrightarrow && \Ac_{1,(f_1)}^\mathsf{lev} \times \Ac_{1,(f_2)}^\mathsf{lev}\\
    (X_1', G_1', G_2') &\mapsto&& (X_1'/G_1', X_1'/G_2')
\end{alignedat}\, .
$$
and $(\varphi_{v,(f_1)} \times \varphi_{g-v, (f_2)}) \circ \Delta$ maps to $\Zc$, and is surjective onto the connected component with discrete data defined above. The reason why $\Delta$ is a "diagonal" map is because if $N= d_1d_2f_1^2f_2^2$ then there is a finite \'etale map
$$
\kappa: \Ac_{1,(N)}^\mathsf{lev} \longrightarrow \Ac_{1,G_1,G_2, f_1, f_2}^{\mathsf{lev}}\, ,
$$
obtained by choosing particular subgroups of $(\Z/N\Z)^2$ which are isomorphic to $G_1$ and $G_2$ and are disjoint, and elements of $(\Z/N\Z)^2$ that give symplectic basis of the $f_1$ and $f_2$-torsion points of the quotients of $(\Z/N\Z)^2$ by these subgroups, and then, we can factor $\Delta \circ \kappa$ as a honest diagonal map
$$
\Ac_{1,(N)}^\mathsf{lev} \longrightarrow \Ac_{1,(N)}^\mathsf{lev} \times \Ac_{1,(N)}^\mathsf{lev}
$$ 
and other finite maps, which quotient by the relevant subgroups.

Going back to the general situation, the normal bundle to the composite morphism
$$
\mathcal W \stackrel{\Delta}{\longrightarrow} \mathcal V \stackrel{\mathrm{mult}_1}{\longrightarrow} \Ac_{v, \mathsf{lev}}\times \Ac_{g-v,\mathsf{lev}}
$$
is
\begin{align*}
    N_{\Delta} + \Delta^* N_{\mathrm{mult}_1} &= \operatorname{Sym}^2 \mathbb E_p^\vee +  
    \mathbb E_{n_2}^\vee \boxtimes \mathbb E_{m_1}^\vee+
    \mathbb E_{n_2}^\vee \boxtimes \mathbb E_{p}^\vee+
    \mathbb E_{p}^\vee \boxtimes \mathbb E_{m_1}^\vee\\
    &+\mathbb E_{n_1}^\vee \boxtimes \mathbb E_{m_2}^\vee+
    \mathbb E_{n_1}^\vee \boxtimes \mathbb E_{p}^\vee+
    \mathbb E_{p}^\vee \boxtimes \mathbb E_{m_2}^\vee\, ,
\end{align*}
And the pullback via $\mathrm{mult}_2 \circ \Delta$ of the normal bundle to
$$
\Ac_{u, \mathsf{lev}}\times \Ac_{g-u, \mathsf{lev}} \longrightarrow\Ac_g
$$
is
\begin{align*}
    (\mathrm{mult}_2 \circ \Delta)^* (\mathbb E_{u}^\vee \boxtimes \mathbb E_{g-u}^\vee)&=(\mathbb E_{n_1}^\vee + \mathbb E_{n_2}^\vee +\mathbb E_{p}^\vee)\boxtimes(\mathbb E_{m_1}^\vee + \mathbb E_{m_2}^\vee +\mathbb E_{p}^\vee)\\
    &=
    \mathbb E_{n_1}^\vee \boxtimes \mathbb E_{m_1}^\vee+
    \mathbb E_{n_1}^\vee \boxtimes \mathbb E_{m_2}^\vee+
    \mathbb E_{n_1}^\vee \boxtimes \mathbb E_{p}^\vee +
    \mathbb E_{n_2}^\vee \boxtimes \mathbb E_{m_1}^\vee\\
    &+\mathbb E_{n_2}^\vee \boxtimes \mathbb E_{m_2}^\vee+
    \mathbb E_{n_2}^\vee \boxtimes \mathbb E_{p}^\vee +
    \mathbb E_{p}^\vee \boxtimes \mathbb E_{m_1}^\vee+
    \mathbb E_{p}^\vee \boxtimes \mathbb E_{m_2}^\vee+
    \mathbb E_{p}^\vee \boxtimes \mathbb E_{p}^\vee\, .
\end{align*}
Therefore, the excess bundle is
$$
\mathrm{Exc}_{\mathcal W} = \mathbb E_{n_1}^\vee \boxtimes \mathbb E_{m_1}^\vee + \mathbb E_{n_2}^\vee \boxtimes \mathbb E_{m_2}^\vee + \wedge^2\mathbb E_p^\vee\, .
$$
Note that if $V,W$ are two vector bundles then
$$
c_{top}(V)=c_{top}(W) = 0 \Rightarrow c_{top}(V \otimes W)=0, \text{ and } c_{top}(V)=0 \text{ or }c_{top}(W)=0 \Rightarrow c_{top}(V \oplus W)=0\, .
$$
The top Chern class of the Hodge bundle vanishes on any moduli space of abelian varieties of dimension $ >0$, so by the Excess intersection formula \cite[Theorem 6.2.]{Ful94}, the Gysin pullback $\varphi_{g,\delta}^!$ in \eqref{eqn: fiber product} is zero when restricted to a component $\mathcal W$ as above with $n_1m_1 \neq 0$ or $n_2m_2 \neq 0$.

Since
$$
\begin{array}{cc}
    u=n_1 +n_2 +p & v=n_2+m_1+p \\
    g-u= m_1 +m_2 +p & g-v = n_1 +m_2 +p
\end{array}\, ,
$$
we cannot have $n_1m_1=n_2m_2=p=0$, so the components $\mathcal W$ where $\varphi_{g,\delta}^!$ could be non-zero satisfy $p\geq 1$.

Let $\mathcal W$ be a component with $n_1m_1 = n_2m_2 =0$, $p >0$. Just as we did in Section \ref{section: 23}, the morphisms $\mathrm{mult}_i, \varphi_{g, \varepsilon}$ can be extended to sufficiently fine toroidal compactifications, while still preserving the splitting of the pullback of the Hodge bundle:
$$
\overline{\mathcal W} \stackrel{\overline{\mathrm{mult}}_i \circ \overline{\Delta}}{\longrightarrow} \overline{\Ac_{v, \mathsf{lev}}\times \Ac_{g-v,\mathsf{lev}}} \stackrel{\overline{\varphi}_{g,\delta}}{\longrightarrow} \overline{\Ac}_g\, .
$$
There is a cycle $\alpha_W$ on $\CH_*(\mathcal{W})$ whose pushforward to $\Ac_g$ is the contribution of $\mathcal{W}$ to $[\NL_{g,\delta}]\cdot [\NL_{g, \varepsilon}]$. If we take an extension $\overline{\alpha}_{\mathcal W}$ of $\alpha_{\mathcal W}$ to $\CH_*(\overline{\mathcal{W}})$, then the sum of the pushforwards of all the $\overline{\alpha}_{\mathcal W}$ will be an extension of $[\NL_{g,\delta}]\cdot [\NL_{g, \varepsilon}]$ to $\CH_*(\overline{\Ac_g})$. The tautological projection of $[\NL_{g,\delta}]\cdot [\NL_{g, \varepsilon}]$ can be computed from integrals of the form
$$
\int_{\overline{\Ac}_g} \sum_{\Wc\text{ with }p \geq 1} (\overline{\mathrm{mult}}_i\circ \overline{\Delta}\circ \overline{\varphi}_{g,\delta})_*(\overline{\alpha}_\Wc)\cdot P(\lambda_i)\cdot c_{top}(\mathbb E_g)\, .
$$
However,
$$
(\overline{\mathrm{mult}}_i\circ \overline{\Delta}\circ \overline{\varphi}_{g,\delta})^*\mathbb E_g = \mathbb E_p^2 \oplus \mathbb E_{n_1} \oplus \mathbb E_{n_2} \oplus \mathbb E_{m_1} \oplus \mathbb E_{m_2} 
$$
has vanishing top Chern class, since $c_p(\mathbb E_p)^2 =0$, on $\overline{\mathcal A}_{p, \mathsf{lev}}$ by Mumford's relation \cite{EV02}.

For arbitrary cycles $\gamma, \eta$ supported on $\NL_{g, \varepsilon}$ and $\NL_{g, \delta}$, $\overline{\gamma \cdot \eta}$ is the pushforward to $\mathcal A_g$ of a cycle supported on components $\overline{\mathcal W}$ with $p>0$, and the same argument applies.
\end{proof}
The intersection product $[\NL_{g, \delta}]\cdot [\NL_{g, \varepsilon}]$ is the sum, over all components $\mathcal W$ of $\mathcal Z$ of a rational multiple of
$$
\mathrm{mult}_{1,*}(\Delta_*(c_{top}( \mathrm{Exc}_{\mathcal W}) \cap \Delta^!(\mathrm{mult}_1^!(1)))) = \begin{cases}
    \mathrm{mult}_{1,*}(\Delta_*(c_{top}(\wedge^2\mathbb E_p^\vee))) &\text{ if }n_1m_1=n_2m_2=0,\\
    0&\text{ otherwise.}
\end{cases}
$$
\begin{lem}\label{lem: diagonal pushforward}
    If $1 \leq p \leq 3$, and
    $$
    \Delta_{\Ac_{p,\mathsf{lev}}} : \Ac_{p,\mathsf{lev}} \longrightarrow\Ac_{p,\mathsf{lev}} \times \Ac_{p,\mathsf{lev}}
    $$
    is the diagonal, then
    $$
    \Delta_{\Ac_{p,\mathsf{lev}},*}(c_{top}(\wedge^2 \mathbb E^\vee_p)) =0\, .
    $$
\end{lem}
\begin{proof}
    The Chow groups of $\Ac_1, \Ac_2, \Ac_3$ have been determined, and they are equal to their tautological rings. Moreover, the stacks $\Ac_1, \Ac_2$ and $\Ac_3$ have the Chow-Kunneth generation property (see \cite[Definition 2.5.]{BS21}) because $\mathcal M_{1,1}$, $\mathcal M_2^{ct}$ and $\mathcal M_3^{ct}$ have it, from their stratifications as quotients of open subsets of affine space, and so
    $$
    \mathsf{R}^*(\Ac_p) \otimes \mathsf{R}^*(\Ac_p) \longrightarrow \CH^*(\Ac_p \times \Ac_p)
    $$
    is surjective. Since $\mathsf{R}^{>{\binom{p}{2}}}(\Ac_p) =0$, we see that $\CH_{<2{\binom{p+1}{2}}-2\binom{p}{2}}(\Ac_p \times \Ac_p)=0$. In particular, $\CH_p(\Ac_p \times \Ac_p)=0$ for $p=1,2,3$. For a moduli stack with level structure, we can form a Cartesian diagram
    $$
        \begin{tikzcd}
	{\Ac_{p,\mathsf{lev}}} & {\Ac_{p,\mathsf{lev}}}\times {\Ac_{p,\mathsf{lev}}} \\
	\Ac_p & {\Ac_p\times \Ac_p}
	\arrow["\pi \times \pi", from=1-2, to=2-2]
	\arrow["\Delta_{\Ac_p}",from=2-1, to=2-2]
	\arrow["\Delta_{\Ac_{p, \mathsf{lev}}}", from=1-1, to=1-2]
	\arrow["\pi", from=1-1, to=2-1]
	\arrow["\lrcorner"{anchor=center, pos=0.125}, draw=none, from=1-1, to=2-2]
\end{tikzcd}.
    $$
    The Hodge bundle pulls back to the Hodge bundle under $\pi$, so the result follows from the formula $(\pi \times \pi)^* \circ \Delta_{\Ac_g,*} = \Delta_{\Ac_{g,\mathsf{lev}}, *} \circ \pi^*$. If $\Ac_{p,\mathsf{lev}}$ does not map to $\Ac_p$, we can always find some $\Ac_{p, \mathsf{lev}'}$ that dominates both, and repeat the argument.
\end{proof}
This proves the following:
\begin{prop}\label{prop: intersection NLd Ag1 Ag2}
    Let $u\leq v \leq g/2$ and $\delta = (d_1, \ldots , d_u)$, $\epsilon_v = (e_1, \ldots , e_v)$ be polarization types. If $u \leq 3$ then
    $$
    [\NL_{g, \delta}]\cdot [\NL_{g, \varepsilon}] =0
    $$
    in the Chow groups of $\Ac_g$.
\end{prop}

\begin{rmk}
    For $p=4,5$, the stacks of Prym curves $\Rc_{p+1}$ have unirational parametrizations, so it is not unreasonable to expect that they have the Chow-Kunneth generation property, and that one can calculate their Chow groups, hoping that they are tautological, and that the same holds for $\Ac_p$ since there are dominant maps $\Rc_{p+1} \longrightarrow\Ac_p$. This would show that Lemma \ref{lem: diagonal pushforward} holds for $p \leq 5$, and so Proposition \ref{prop: intersection NLd Ag1 Ag2} is likely true for $u \leq 5$. It would be interesting to see if this is true for any $u$, and this requires a new idea, since $\Ac_6$ is not unirational \cite{DMS21} and has non-tautological algebraic classes \cite{COP}.
\end{rmk}

\section{The Hecke algebra}\label{section: 4}

For each $n \geq 1$, let
$$
\mathcal B_{g, n} = \left\lbrace \begin{array}{c}
     \text{homomorphisms }f:(X, \theta_X) \longrightarrow(Y, \theta_Y)\\
     \text{ such that } (X, \theta_X), (Y, \theta_Y) \in \Ac_g\text{ and }f^*\theta_Y = n\theta_X 
\end{array}\right\rbrace.
$$
It is a smooth DM stack and has forgetful maps
$$
\begin{tikzcd}
	& \Bc_{g,n} \\
	{\Ac_g} && {\Ac_g}
	\arrow["\pi_1"', from=1-2, to=2-1]
	\arrow["\pi_2", from=1-2, to=2-3]
\end{tikzcd}
$$
which are \'etale, and induce \emph{Hecke operators}
$$
T_n := \pi_{2,*} \circ \pi_1^* :  \CH^*(\Ac_g)\longrightarrow\CH^*(\Ac_g)\, .
$$
Let $\vol(T_n) = T_n(\mathbf{1}) \in \Z$ be the \textit{volume} of $T_n$. We compute in the Appendix (Proposition \ref{prop:volumeTn}) that
$$
\vol (T_n) = \prod_{p \mid n}\prod_{i=1}^g\frac{1-p^{g(v_p(n)+i)}}{1-p^{gi}}\, .  
$$
It is clear that $T_n$ maps tautological classes to tautological classes, since the Hodge bundle is invariant under isogenies. In fact, they are eigenvectors of all the $T_n$ with eigenvalues $\vol(T_n)$.

\begin{lem}\label{lem: taut commutes hecke}
    The tautological projection commutes with the Hecke operators.
\end{lem}
\begin{proof}
    $\pi_1$ and $\pi_2$ are morphisms of locally symmetric domains, so by what we explained in Section \ref{section: 23}, they extend to sufficiently fine toroidal compactifications, and the extensions satisfy
$\overline{\pi}_1^*\mathbb E_g = \overline{\pi}_2^* \mathbb E_g$, so if $\alpha \in \CH^*(\Ac_g)$ and $\Lambda$ is a tautological class,
\begin{align*}
    \int_{\overline{\Ac}_g} \taut\left( T_n(\alpha)\right)\Lambda \lambda_g &= \int_{\overline{\Ac}_g}\overline{T_n(\alpha)} \Lambda \lambda_g = \int_{\overline{\Ac}_g}T_n(\overline{\alpha}) \Lambda \lambda_g =\int_{\overline{\mathcal A}_g} \overline{\pi}_{2,*}(\overline{\pi}_1^*\alpha)\Lambda \lambda_g\\
    &= \int_{\overline{\mathcal B}_{g,n}} \overline{\pi}_1^*\alpha \cdot \overline{\pi}_{2}^*(\Lambda \lambda_g) = \int_{\overline{\mathcal B}_{g,n}} \overline{\pi}_1^*(\alpha\cdot \Lambda \lambda_g)=\vol(T_n) \int_{\overline{\Ac}_g}\overline{\alpha} \Lambda \lambda_g \\
    &= \vol(T_n) \int_{\overline{\Ac}_g}\taut(\alpha) \Lambda \lambda_g = \int_{\overline{\Ac}_g}T_n(\taut(\alpha)) \Lambda \lambda_g\, .
\end{align*}
\end{proof}
Abelian subvarieties (and their dimensions) are preserved by isogenies, so Hecke operators preserve the $\Q$-span of these $\NL$-cycles. We focus on the $\NL$-cycles corresponding to elliptic subgroups.

\begin{defn}
    Let
    $$
    \mathbf{NL}_g \subset \CH^{g-1}(\Ac_g)
    $$
    be the subspace generated by $\lambda_{g-1}$ and all the $\NL_{g,d}$ for $d \geq 1$. The NL-\textit{Hecke algebra} is the subalgebra
    $$
        \mathsf{Hecke}^{\NL}_g  \subset \End_\Q(\mathbf{NL}_g)
    $$
    generated by all the $T_n$.
\end{defn}

\begin{prop}\label{thm: Hecke algebra}
    The NL-Hecke algebra has the following properties:
    \begin{enumerate}[label = \alph*)]
        \item $\mathsf{Hecke}^{\NL}_g$ is commutative, preserves the $\Q$-span of $\lambda_{g-1}$.
        \item As a $\mathsf{Hecke}^{\NL}_g$-module, $\mathbf{NL}_g$ is generated by $[\Ac_1 \times \Ac_{g-1}]$ and $\lambda_{g-1}$.
        \item $\mathsf{Hecke}^{\NL}_g$ is is generated by the $T_p$ for $p$ prime.
    \end{enumerate}
    In particular, if $[\Ac_1\times \Ac_{g-1}]$ is tautological, then $\mathbf{NL}_g$ is $1$-dimensional.
\end{prop}
\begin{proof}
    The Hecke algebra of correspondences on $\Ac_g$ is well-known to be commutative, see \cite[Chapter VII]{CF91}, and the rest of a) has been discussed before. For b), note that
    $$
    \pi_2\left(\pi_1^{-1}(\Ac_1 \times \Ac_{g-1})\right)
    $$
    is the set of principally polarized abelian varieties $(Y, \theta_Y)$ having an isogeny $f: E \times X' \longrightarrow Y$, where $X'$ is also principally polarized and $f^*\theta_Y = n\theta_E + n\theta_{X'}$. The induced map
    $$
    g: E \longrightarrow E \times \{0\} \subset E \times X' \longrightarrow Y
    $$
    satisfies that $g^*\theta_Y = n \theta_E$, and so $g(E)$ is an elliptic subgroup of $Y$ and $\theta_{Y | g(E)}$ has degree $n/|\ker (g)|$.
    
    Given $Y$ with an elliptic subgroup $F$ such that $\theta_{Y|F}$ has degree $m \mid n$, let $Y'$ be the complementary subvariety of $F$. $\theta_{Y|Y'}$ is of type $(1,\ldots , 1,m)$, so in particular, there is a principally polarized $X'$ and an isogeny $X' \longrightarrow Y'$ such that the pullback of $\theta_{Y'}$ is $n$ times the polarization of $X'$. If $E \longrightarrow F$ has degree $n/m$ then $X:= E \times X'$ maps to $Y$ under the correspondence. Therefore,
    $$
    \pi_2\left(\pi_1^{-1}(\Ac_1 \times \Ac_{g-1})\right) = \bigcup_{m \mid n}\NL_{g,m}\, ,
    $$
    and 
    $$
    T_n([\Ac_1 \times \Ac_{g-1}]) = \sum_{m \mid n}c_{n,m} [\NL_{g,m}]\, .
    $$
    for some positive coefficients $c_{n,m}$. This proves b).

    Any isogeny decomposes uniquely into isogenies of relatively prime degrees. Therefore, $\mathsf{Hecke}_g^{\NL}$ is generated by all the $T_{p^k}$. By the same arguments that we explained before,
    $$
    \pi_2\left(\pi_1^{-1}(\NL_{g,d})\right) = \bigcup_{m \mid nd}\NL_{g,m}
    $$
    and so,
    $$
    T_p(\stackrel{l \text{ times}}{\ldots} (T_p([\Ac_1 \times \Ac_{g-1}])\ldots ) = \sum_{i=0}^l d_{l,i} [\NL_{g,p^i}]
    $$
    and 
    $$
    T_{p^k}([\Ac_1 \times \Ac_{g-1}]) = \sum_{i=0}^k e_{k,i} [\NL_{g,p^i}]
    $$
    for some positive numbers $d_{l,i}, e_{k,i}$.
    Therefore, one solves a triangular system to find a polynomial $q \in \Q[T]$ of degree $k$ such that
    $$
    (T_{p^k} - q(T_p)) ([\Ac_1 \times \Ac_{g-1}])=0\, .
    $$
    By Lemma \ref{lem: taut commutes hecke} the Hecke operators commute with $\taut$, and $\taut([\Ac_1 \times \Ac_{g-1}])$ is a multiple of $\lambda_{g-1}$. Therefore,
    $$
    (T_{p^k} - q(T_p)) (\lambda_{g-1}) =0\, .
    $$
    Since the Hecke algebra is commutative, and $\mathbf{NL}_g$ is generated by $\lambda_{g-1}$ and $[\Ac_1 \times \Ac_{g-1}]$, this shows that the relation $T_{p^k} = q(T_p)$ holds in $\mathsf{Hecke}_g^{\NL}$.

    If $[\Ac_1\times \Ac_{g-1}]$ is tautological, it must equal its tautological projection, so $\mathbf{NL}_g$ would be generated by $\lambda_{g-1}$ and the Hecke operators, but $\lambda_{g-1}$ is a common eigenvector of them.
    \end{proof}
The space $\operatorname{Mod}_{2g}$ of modular forms of weight $2g$ for $\operatorname{SL}_{2}(\mathbb Z)$ also has an action of a commutative Hecke algebra generated by elements $T_p$ for primes $p$, and is generated by an Eisenstein form and a cusp form as a module over the Hecke algebra. Therefore, Proposition \ref{thm: Hecke algebra} supports Conjecture \ref{conj: modularity}.

\begin{proof}[Proof of Theorem \ref{thm: not tautological}]
    Note that $\mathrm{Tor}^{-1}(\NL_{g,2}) \cap \mathcal M_g$ is precisely the bielliptic locus, which has the expected codimension $g-1$, so the image of $[\NL_{g,2}]$ by the maps
    $$
    \CH^{g-1}(\Ac_{g}) \stackrel{\Tor^*}{\longrightarrow} \CH^{g-1}(\Mc_g^{ct}) \longrightarrow \CH^{g-1}(\Mc_g)
    $$
    is a positive multiple of the class of the bielliptic locus, which is not tautological by the work of \cite{ACHLMT24} in the cases $g=12$ and $g \geq 16$ even. However, these maps send tautological classes to tautological classes. Therefore, $\NL_{g,2}$ cannot be tautological and by Proposition \ref{thm: Hecke algebra}, $[\Ac_1\times \Ac_{g-1}]$ cannot be tautological either.
\end{proof}

\begin{rmk}
In fact, since $\Tor^{-1}(\Ac_1 \times \Ac_{g-1}) \cap \Mc_g = \emptyset$, we cannot write $[\NL_{g,2}]$ in terms of $[\Ac_1 \times \Ac_{g-1}]$ and $\lambda_{g-1}$. Therefore, $\dim(\mathbf{NL}_g) \geq 3$ for $g=12$ or $g \geq 16$ even.
\end{rmk}

\newpage
\section{The homomorphism property}\label{section: 5}

Recall from Definition \ref{defn: homomorphism property} that $(\alpha, \beta)$ has the homomorphism property if
$$
\taut(\alpha \cdot \beta) = \taut(\alpha)\cdot \taut(\beta)\, .
$$

\begin{proof}[Proof of Theorem \ref{thm: homomorphism property}]
Note that
$$
\lambda_{g-1}\mid_{\NL_{g, \delta}} =0\, .
$$
this is because, after pulling back to $\Ac_{u, g-u, \delta}$, $\mathbb E_g$ splits as the sum of two vector bundles with vanishing top chern class. In \cite{vdG99} it is explained that $\mathsf{R}^*(\Ac_g)$ is additively generated over $\Q$ by monomials $\lambda_1 ^{u_1}\ldots \lambda_{g-1}^{u_{g-1}}$ with $u_i \in \{0,1\}$, so if $\alpha \in \CH^*(\Ac_g)$ is a class such that $\lambda_{g-1}\cdot \alpha = 0$ then $\lambda_{g-1}\cdot \taut(\alpha) = 0$ and so $\lambda_{g-1} \mid \taut(\alpha)$.

This shows that $\lambda_{g-1} \mid \taut(\alpha)$ whenever $\alpha$ is supported on a $\NL$ cycle, so $ 0 = \lambda_{g-1}^2 \mid \taut(\alpha) \cdot \taut(\beta)$ whenever $\alpha$ and $\beta$ are supported on $\NL$ cycles. Finally, $\taut(\alpha \cdot \beta)$ was shown to be $0$ in Proposition \ref{prop: taut of a product}.
\end{proof}

The homomorphism property is connected to the $\lambda_g$-pairing on $\mathcal{M}_g^{ct}$ by the following, which is a generalization of the proof of \cite[Theorem 5]{COP}.
\begin{prop}\label{prop: homomorphism and Gorenstein}
    Let $k \leq 2g-3$ be a positive integer such that $R^{2g-3-k}(\Mc_g)$ is generated by monomials in the $\lambda_i$\footnote{For example, $k=2g-3, 2g-2$ or any $k \leq g-1$, or any $k$ for small values of $g$.}. If $\alpha \in \CH^k(\Ac_g)$ and $(\alpha , \beta)$ has the homomorphism property whenever $\beta$ is supported on $\mathcal A_u \times \mathcal A_{g-u}$ for all $u$, then the pair $(\alpha, \Tor_*(\mathbf{1}))$ has the homomorphism property if and only if $\Tor^*(\alpha - \taut(\alpha)) \in \CH^k(\Mc_g^{ct})$ lies in the kernel of the $\lambda_g$-pairing
    $$
    \CH^k(\Mc_g^{ct}) \times \mathsf{R}^{2g-3-k}(\Mc_g^{ct}) \longrightarrow \CH^{2g-3}(\Mc_g^{ct}) \longrightarrow \Q\, .
    $$
\end{prop}
\begin{proof}
    We work on a toroidal compactification of $\Ac_g$ such that the Torelli map extends to a proper morphism $\overline{\Tor}$ defined on $\overline{\Mc}_g$; for example, the second Voronoi compactification works by \cite{Nam76}, and the  perfect cone compactification also works by \cite{AB12}. By formula \eqref{eqn: semihomomorphism},
    $$
    \taut(\alpha \cdot \Tor_*(\mathbf{1})) - \taut(\alpha) \cdot \taut(\Tor_*(\mathbf{1})) = \taut(\Tor_*(\mathbf{1})\cdot (\alpha -\taut(\alpha))\, .
    $$
    so this class is zero if and only if
    $$
    \int_{\overline{\Mc}_g} (\overline{\Tor}^*(\overline{\alpha}) - \overline{\Tor}^*(\taut(\alpha)))\Lambda \lambda_g=0\, .
    $$
    for any polynomial $\Lambda$ on the $\lambda_i$. We are using that $\overline{\Tor}^*(\overline{\alpha})$ is an extension of $\Tor^*(\alpha)$ to $\overline{\Mc}_g$ and that $\overline{\Tor}^*$ preserves the Hodge bundle.
    
    By the assumption on $k$, the rest of the classes in $\mathsf{R}^{2g-3-k}(\Mc_{g}^{ct})$ are supported on the boundary $\Mc^{ct}_g \smallsetminus \Mc_g$, but their pushforward under the Torelli map lies in $\mathcal A_{u}\times \mathcal A_{g-u}$ for some $u$, so the homomorphism property holds for these pushforwards by Theorem \ref{thm: homomorphism property}.
\end{proof}

\section{Gromov-Witten theory of the universal elliptic curve}\label{section: 6}

\begin{defn}\label{defn: relative stable maps}
    Let $\pi : \overline{\mathcal E} \longrightarrow\overline{\Mc}_{1,1}$ be the universal curve, with section $q$. The $\pi$-relative moduli space of log-stable maps of degree $d$ from genus $g$ curves is denoted by
    $$
    \overline{\Mc}_{g,1}(\pi, d)\, .
    $$
    Over a scheme $S$, it parametrizes the data of:
    \begin{itemize}
        \item a nodal curve $C \longrightarrow S$ of genus $g$, with its canonical log structure (given by its singular locus),
        \item a quasi-stable\footnote{A pointed nodal curve is quasi-stable if it has no unstable rational tails} curve $E \longrightarrow S$ of arithmetic genus $1$ with a section, and the canonical log structure given by its singular locus,
        \item a log stable map $f:C \longrightarrow E$ such that $\operatorname{Aut}(f)$ is finite.
    \end{itemize}
\end{defn}
The log condition says that, if $E$ is a singular elliptic curve, then the preimage by $f$ of a node consists of nodes, and on an \'etale neighborhood of each node in the domain, the map is given by
$$
\begin{array}{rcl}
    \frac{\C[x,y]}{xy=0} &\longrightarrow & \frac{\C[a,b]}{ab=0}  \\
    (x,y) & \mapsto & (a^k, b^k)\, .
\end{array}
$$
This space comes with a virtual fundamental class of dimension
$$
\mathrm{vdim}\left(\overline{\Mc}_{g,1}(E,d)\right) + \dim\left(\overline{\Mc}_{1,1}\right) =2g\, ,
$$
and natural morphisms
$$
\mathsf{ev}:\overline{\Mc}_{g,1}(\pi, d) \longrightarrow \overline{\mathcal E}\text{ and }\mathsf{ft} :\overline{\Mc}_{g,1}(\pi, d) \longrightarrow \overline{\Mc}_{g,1}\, .
$$
Let $\Mc_{g,1}^{ct}(\pi,d)$ be the inverse image of $\Mc^{ct}_{g,1}$ under $\mathsf{ft}$, and $\mathcal E$ be the universal elliptic curve over $\Mc_{1,1}$.
\begin{lem}\label{lem: follow the cycle}
    If $d >0$, the map $\mathsf{ev} : \Mc_{g,1}^{ct} (\pi,d) \longrightarrow\overline{\Ec}$ factors through $\Ec$.
\end{lem}
\begin{proof}
    Let $P = P_1 \cup \ldots \cup P_k$ be a cycle of $k$ rational curves, with nodes $n_1, \ldots , n_k$. If $f:C \longrightarrow P$ is a log-stable map and $f(q) =n_i$ then $q$ has to be a node. If $C' \subset C$ is an irreducible component that is not contracted by $f$, we obtain a surjective map from $C'$ to $P_i$ for some $i$, so there are at least two points in $C'$ mapped to some node. Therefore, $f$ contracts all the irreducible components of $C$ whose corresponding vertex in the dual graph is adjacent to only one edge. If $C$ is of compact-type, one can contract all the leaves of its dual graph repeatedly, and so $f$ has to be constant, which is impossible because $d>0$.
\end{proof}

\begin{lem}\label{lem: psi class}
    $\mathsf{ft}^*(\psi_1)$ agrees with the $\psi_1$ class on $\mathcal{M}^{ct}_{g,1}(\pi, d)$.
\end{lem}
\begin{proof}
    Indeed, by the Lemma \ref{lem: follow the cycle}, every rational component of the domain of a stable map in $\Mc_{g,1}^{ct}(\pi, d)$ maps to an elliptic curve, so it is contracted. Therefore the domain curve is always stable, so the two $\psi$ classes are the same.
\end{proof}
Note that for a pointed compact type curve $(C, p)$ and an elliptic curve,
\begin{equation}\label{eqn: correspondence NL maps}
\left\lbrace \begin{array}{cc}
   \text{Homomorphisms }f : E \longrightarrow\operatorname{Jac}(C)    \\
     \text{such that}\deg (f^*(\theta_C))=d
\end{array}\right\rbrace \longleftrightarrow \left\lbrace \begin{array}{cc}
   \text{Maps }h : (C,p) \longrightarrow(E,0)    \\
     \text{of degree }d
\end{array}\right\rbrace. 
\end{equation}
Indeed, given $f:E \longrightarrow\operatorname{Jac}(C)$ and $p$, the map $h$ is $f^\vee \circ \operatorname{aj}_p$, where
$$
\operatorname{aj}_p : (C,p) \longrightarrow (\operatorname{Jac}(C), 0)
$$
is the Abel-Jacobi map, and given $h$, $f$ will be the dual of the map $\operatorname{Jac}(C) \longrightarrow E$ given by the Albanese property of $\operatorname{Jac}(C)$.

We construct the moduli spaces of both sides of the equivalence \eqref{eqn: correspondence NL maps}. Let $\mathcal{Q}_{g,d}$ be the fibered product
$$
\begin{tikzcd}
	\mathcal{Q}_{g,d} & {\Mc^{ct}_{g,1}(\pi, d)} \\
	{\Mc_{1,1}} & \Ec
	\arrow[from=1-1, to=1-2]
	\arrow[from=1-1, to=2-1]
	\arrow["\lrcorner"{anchor=center, pos=0.125}, draw=none, from=1-1, to=2-2]
	\arrow["{\mathsf{ev}}", from=1-2, to=2-2]
	\arrow["q"', from=2-1, to=2-2]
\end{tikzcd}.
$$
Consider the morphism
$$
\Tor_1 \colon \Mc_{g,1}^{ct} \longrightarrow \Ac_g\, ,
$$
which is the composition of the Torelli morphism with the forgetful morphism. We define $\Tor_1^{-1}(\widetilde{\NL}_{g,d})$ as the fibered product
$$
\begin{tikzcd}
	{\Tor_1^{-1}(\widetilde{\NL}_{g,d})} & {\Mc^{ct}_{g,1}} \\
	{\widetilde{\NL}_{g,d}} & {\Ac_g}
	\arrow[from=1-1, to=1-2]
	\arrow[from=1-1, to=2-1]
	\arrow["\lrcorner"{anchor=center, pos=0.125}, draw=none, from=1-1, to=2-2]
	\arrow["{\Tor_1}", from=1-2, to=2-2]
	\arrow[from=2-1, to=2-2]
\end{tikzcd}
$$
where $\widetilde{\NL}_{g,d}$ is the stack of maps $E \longrightarrow X$ of degree $d$ (and was already considered in \ref{cor: modular form}). 
By \eqref{eqn: correspondence NL maps}, which clearly works in families, the following holds:

\begin{lem}[Generalizing \cite{COP}, Proposition 21, see also \cite{GL24}]\label{prop: Tor to stable 1}
    There is an isomorphism of stacks
    $$
    \Qc_{g,d} \longrightarrow \Tor_1^{-1}(\widetilde{\NL}_{g,d})\, .
    $$
    over $\Mc_{g,1}^{ct}$.
\end{lem}
\begin{proof}
    The morphism $\Qc_{g,d} \longrightarrow\widetilde{\NL}_{g,d}$ is the one that we have described above using the Abel-Jacobi map. Together with $\mathsf{ft}$, this gives a natural morphism $\Qc_{g,d} \longrightarrow \Tor_1^{-1}(\widetilde{\NL}_{g,d})$, and an inverse is constructed likewise, because the equivalence \eqref{eqn: correspondence NL maps} works in families.
\end{proof}
But moreover, there is an equality of virtual cycles:
\begin{lem}[See \cite{GL24}, Theorem 2]\label{thm: Tor to stable 2}
    Under the isomorphism in Lemma \ref{prop: Tor to stable 1}, there is an equality of cycles
    $$
    \Tor_1^*([\widetilde{\NL}_d]) = q^!\left(\left[\Mc_{g,1}^{ct}(\pi,d)\right]^{\mathrm{vir}}\right)\, .
    $$
    In particular, there is a canonical extension of $\Tor_1^*([\widetilde{\NL}_{g,d}])$ to $\overline{\Mc}_{g,1}$.
\end{lem}

\begin{proof}[Proof of Theorem \ref{thm: equivalent conj}]
    First note that $[\NL_{g,d}]$ and $[\widetilde{\NL}_{g,d}]$ are related by the invertible relation \eqref{eqn: tilde to not tilde}, so we may exchange them freely. Recall also that the tautological projection of $\widetilde{\NL}_{g,d}$ was calculated in the proof of Corollary \ref{cor: modular form}.
    
    Since $[\widetilde{\NL}_{g,d}]$ is a class of codimension $g-1$, has the homomorphism property with respect to any class supported on a $\NL$-cycle (\ref{thm: homomorphism property}), and $\mathsf{R}^{g-2}(\Mc_g)$ is generated by $\lambda$ classes, we may apply \ref{prop: homomorphism and Gorenstein}, so $[\widetilde{\NL}_{g,d}]$ and $\Tor_*(1)$ have the homomorphism property if and only if
    \begin{equation}\label{eqn: tor pullback}
    \Tor^*\left([\widetilde{\NL}_{g,d}] - \frac{g\sigma_{2g-1}(d)}{6|B_{2g}|}\lambda_{g-1} \right)
    \end{equation}
    pairs to $0$ with every tautological class on $\Mc_g^{ct}$. By Lemmas \ref{thm: Tor to stable 2} and \ref{lem: psi class},
    $$
    \frac{1}{2g-2}\pi_*\mathsf{ft}_*\left( q^!\psi_1\left[\overline{\Mc}_{g,1}(\pi,d)\right]^{\mathrm{vir}}\right)
    $$
    is an extension of $\Tor^*([\widetilde{\NL}_{g,d}])$ to $\overline{\Mc}_g$, where the $\frac{1}{2g-2}$ factor comes from the insertion of the $\psi$ class. In fact, $\mathsf{R}^{g-2}(\Mc_g)$ is generated by $\lambda_{g-2}$ by the work in \cite{Loo95}, so \eqref{eqn: tor pullback} is in the kernel of the $\lambda_g$ pairing if and only if
    $$
    \frac{1}{2g-2} \int_{[\overline{\Mc}_{g,1}(\pi,d)]^{\mathrm{vir}}} \psi_1 \mathsf{ev}^*(q) \lambda_{g-2}\lambda_g = \frac{g\sigma_{2g-1}(d)}{6|B_{2g}|} \int_{\overline{\Mc}_g}\lambda_{g-2}\lambda_{g-1}\lambda_g = \frac{|B_{2g-2}|\sigma_{2g-1}(d)}{24(2g-2)(2g-2)!}\, ,
    $$
    where the last integral was computed in \cite{FP00}. This completes the proof.
\end{proof}

\newpage
\section{Degree computations}\label{section: 7}

We compute the degrees of the maps $\phi_\delta$ and $\pi_\delta$ and the volume of the Hecke operator $T_n$.

\subsection{The degree of $\phi_\delta$}\label{section: 71}

We first consider the case $\delta = (1,\ldots ,1, d, \ldots ,d)$ to illustrate the ideas that we will use.
\subsubsection{$\delta = (1, \ldots , 1,d,\ldots ,d)$}

Throughout this part, $k$ is the number of $1$ entries and $h$ is the number of $d$ entries, so $h +k=g$.
\begin{prop}\label{prop: degree of aD}
    When $\delta = (\underbrace{1, \ldots ,1}_{k \text{ times}},\underbrace{d, \ldots ,d}_{h \text{ times}})$, the degree of $\phi_\delta$ is
    $$
    d^{h(2g+1)} \prod_{p \mid d}\prod_{i=g-h+1}^{g}(1-p^{-2i})\, ,
    $$
    where the product is over all primes dividing $d$.
\end{prop}

By Lemma \ref{lem: alpha over the complex}, this is the same as the index of the subgroup of $\Sp_{2g}(\Z)$ consisting of matrices that can be written as
\begin{equation}\label{eqn: form matrices}
\mathbbm{1}_{2g} +  \begin{bmatrix}
a_{11} & \ldots & a_{1g} & b_{11} & \ldots & b_{1k} & db_{1,k+1} & \ldots & db_{1g} \\
\vdots & \ddots & \vdots & \vdots & \ddots & \vdots & \vdots & \ddots & \vdots \\
a_{k1} & \ldots & a_{kg} & b_{k1} & \ldots & b_{kk} & db_{k, k+1} & \ldots & db_{kg} \\
da_{k+1,1} & \ldots & da_{k+1,g} & db_{k+1,1} & \ldots & db_{k+1,k} & d^2b_{k+1,k+1} & \ldots & d^2b_{k+1,g} \\
\vdots & \ddots & \vdots & \vdots & \ddots & \vdots & \vdots & \vdots & \vdots \\
da_{g1} & \ldots & da_{gg} & db_{g1} & \ldots & db_{gk} & d^2b_{g,k+1} & \ldots & d^2b_{gg} \\
c_{11} & \ldots & c_{1g} & e_{11} & \ldots & e_{1k} & de_{1,k+1} & \ldots & de_{1g} \\
\vdots & \ddots & \vdots & \vdots & \ddots & \vdots & \ldots & \ddots & \vdots \\
c_{g1} & \ldots & c_{gg} & e_{g1} & \ldots & e_{gk} & de_{g,k+1} & \ldots & de_{gg} 
\end{bmatrix}\, ,
\end{equation}
where the $a_{ij}, b_{ij}, c_{ij}, e_{ij}$ are integers.

For an integer $n$ and a prime $p$, $v_p(n)$ denotes the highest power of $p$ that divides $n$. If we write $d = p_1^{v_{p_1}(d)}\cdot \ldots \cdot p_m^{v_{p_m}(d)}$,
and $\delta_i = (1,\ldots ,1, p_i^{v_{p_i}(d)}, \ldots , p_i^{v_{p_i}(d)})$ then
\begin{equation}\label{eqn: localize at primes}
\Gamma_{g}[\delta] = \bigcap_{i=1}^m \Gamma_{g}[\delta_i]\, ,
\end{equation}
and so we can reduce to the case where $d=p^v$.

We are going to stratify $\Grm_\delta[\delta]$ by intersecting it with the normal subgroups
$$
\Gamma_g[p^i] = \Sp_{2g}(\Z)[p^i] = \ker \left(\Sp_{2g}( \Z) \longrightarrow\Sp_{2g}(\Z/p^i\Z)\right)\, ,
$$
and since $\Gamma_g[d^{2v}] \leq \Grm_\delta[\delta]$,
$$
[\Sp_{2g} (\Z) : \Grm_\delta[\delta]] \cdot \prod_{i=1}^{2v}\left|\mathrm{im}\left(\Grm_\delta[\delta] \cap \Gamma_g[p^{i-1}] \longrightarrow\Sp_{2g}( \mathbb Z/p^i\Z)\right)\right| = \left|\Sp_{2g}( \Z/p^{2v}\Z) \right|\, .
$$
After stratifying $\Gamma_{g}[p^{2v}]$ in the same way, one arrives at the equality
\begin{equation}\label{eqn: general group quotient}
    [\Sp_{2g}(\Z) : \Grm_\delta[\delta]] =\prod_{i=1}^{2v}\frac{\left|\mathrm{im}\left(\Gamma_{g}[p^{i-1}] \longrightarrow\Sp_{2g}(\Z/p^{i}\Z)\right)\right|}{\left|\mathrm{im}\left(\Grm_\delta[\delta] \cap \Gamma_g[p^{i-1}] \longrightarrow\Sp_{2g}(\mathbb Z/p^i\Z)\right)\right|}
\end{equation}
The $i=0$ computation is different, so we do this first.
By a Gram-Schmidt procedure, any linearly independent set of vectors that spans an isotropic subspace of $\mathbb F_p^{2g}$ can be extended to a symplectic basis. In other words, $\Sp_{2g}(\mathbb F_p)$ acts transitively on the set
$$
\{(e_1, \ldots , e_h) \in \mathbb (\mathbb F_p^{2g})^h: \langle e_1, \ldots , e_h\rangle \text{ is an }h\text{-dimensional isotropic subspace}\}\, ,
$$
which has cardinality
$$
(p^{2g}-1)(p^{2g-1}-p)\ldots (p^{2g-(h-1)}-p^{h-1}) = p^{2gh-\binom{h}{2}} \prod_{i=g-h+1}^g(1-p^{-2i})
$$
From \eqref{eqn: form matrices}, we see that
$$
\mathrm{im}\left(\Grm_\delta[\delta] \longrightarrow\Sp_{2g}(\mathbb F_p)\right)
$$
is the same as the set of symplectic endomorphisms $\mathbb F_p^{2g} \longrightarrow\mathbb F_p^{2g}$ that leave fixed the last $h$ vectors of the standard symplectic basis; in other words, it is the stabilizer of a point for the above action, and so
$$
\frac{|\Sp_{2g}( \mathbb F_p)|}{\left|\mathrm{im}\left(\Grm_\delta[\delta] \longrightarrow\operatorname{Sp}_{2g}( \mathbb F_p)\right)\right|} = p^{2gh-\binom{h}{2}}\prod_{i=g-h+1}^g(1-p^{-2i})\, .
$$

For $i>1$, an element of $\mathrm{im}\left(\Gamma_g[p^{i-1}] \longrightarrow\Sp_{2g}( \Z/p^i\Z)\right)$ can be uniquely written as
\begin{equation}\label{eqn: induction step}
    M=\mathbbm{1}_{2g} + p^{i-1}\begin{pmatrix}
    A&B\\C&-A^t
\end{pmatrix}\, ,
\end{equation}
where $A,B,C,D$ have entries in $\mathbb F_p$, and $B,C$ are symmetric. If we write
$$
A = \begin{pmatrix}
    A_1\\
    A_2
\end{pmatrix}\, , \quad B = \begin{pmatrix}
    B_1&B_2\\
    B_2^t&B_4
\end{pmatrix}\, ,
$$
where $A_1\in \mathrm{M}_{k\times g}(\mathbb F_p)$, $A_2\in \mathrm{M}_{h\times g}(\mathbb F_p)$, $B_1\in \operatorname{Sym}_{k\times k}(\mathbb F_p)$, $B_2\in \mathrm{M}_{k\times h}(\mathbb F_p)$, $B_4\in \operatorname{Sym}_{h\times h}(\mathbb F_p)$. After inspecting \eqref{eqn: form matrices} we arrive at the following:
\begin{itemize}
    \item If $1<i\leq v$ then $M$ is in $\mathrm{im}\left(\Grm_\delta[\delta]\cap\Gamma_g[p^{i-1}] \longrightarrow\Sp(2g, \Z/p^i)\right)$ if and only if $A_2=B_2=B_4=0$ over $\mathbb F_p$.
    \item If $v<i\leq 2v$ then $M$ is in $\mathrm{im}\left(\Grm_\delta[\delta]\cap\Gamma_g[p^{i-1}] \longrightarrow\Sp(2g, \Z/p^i)\right)$ if and only if $B_4=0$ over $\mathbb F_p$.
\end{itemize}
Therefore
$$
\frac{\left|\mathrm{im}\left(\Gamma_{g}[p^{i-1}] \longrightarrow\Sp(2g, \Z/p^{i}\Z)\right)\right|}{\left|\mathrm{im}\left(\Grm_\delta[\delta] \cap \Gamma_g[p^{i-1}] \longrightarrow\operatorname{Sp}(2g, \mathbb Z/p^i\Z)\right)\right|} = \begin{cases}
    p^{hg + hk + h(h+1)/2} &\text{if }1<i\leq v\\
    p^{h(h+1)/2}&\text{if }v<i\leq 2v
\end{cases}\, ,
$$
and this proves Proposition \ref{prop: degree of aD}, since the sum of all the powers of $p$ that appear is
$$
p^{2gh-\binom{h}{2} + (v-1)(hg + hk + h(h+1)/2)  + (2v-v)(h(h+1)/2)} = p^{vh(2g+1)}\, .
$$

\subsubsection{The general case}

For a general sequence $\delta$, by \eqref{eqn: localize at primes} we can assume that $\delta =(1,\ldots , 1, p^{v_1}, \ldots , p^{v_h})$ is formed by powers of a fixed prime $p$, with $1\leq v_1\leq \ldots \leq v_h$. The formula \eqref{eqn: general group quotient} is unchanged, with $l=v_h$. The $i=1$ factor of the right hand side of \eqref{eqn: general group quotient} is again
$$
p^{2gh-\binom{h}{2}}\prod_{i=g-h+1}^g(1-p^{-2i})\, ,
$$
and the rest of the factors can be are obtained by stratifying. Note that on the $i$-th step of the computation of the stratification, we are looking at matrices of the form \eqref{eqn: induction step}, its contribution to the product \eqref{eqn: general group quotient} is $p^{N(i)}$ where $N(i)$ is the number of entries $a_{ij}'$ (for any $i,j$) or $b_{ij}'$ (for $i \leq j$) of a matrix 
$$
M' = \left(\begin{array}{c|c}
    (a_{ij}') & (b_{ij}') \\
    \hline
    (c_{ij}') & (e_{ij}')
\end{array}\right)
$$
of the form \eqref{eqn: form matrices} (with $\delta = (1,\ldots , 1, p^{v_1}, \ldots , p^{v_h})$) that are always divisible by $p^i$. A simple computation shows that
$$
\sum_{i=1}^{2v_h} N(i) = (v_1 + v_2+\ldots + v_h) (2g-h) + \sum_{i \leq j} (v_i +v_j) = (2g+1)(v_1 + \ldots + v_h),
$$
and $N(1) = h(2g-h) + h(h+1)/2 = 2gh - h(h-1)/2$. Therefore,
\begin{align*}
    [\Sp_{2g}(\Z) : \Grm_\delta[\delta]] &= \left(p^{2gh-\binom{h}{2}}\prod_{i =g-h+1}^g(1-p^{-2i})\right) p ^{N(2) + \ldots + N(2v_h)}\\
    &= p^{(2g+1)(\sum_j v_j) } \prod_{i=g-h+1}^g(1-p^{-2i})\, .
\end{align*}

This shows that:
\begin{prop}\label{prop: degree aD}
    If $\delta = (d_1, \ldots , d_g)$ is a polarization type with $d=d_1\ldots d_g$, denote by
    $$
    h_p(\delta) = g-\max\{i : v_p(d_i)=0\}\, .
    $$
    Then, the degree of $\phi_\delta$ is 
    $$
    d^{2g+1} \prod_{p \mid d} \prod_{i=g-h_p(\delta)+1}^{g}(1-p^{-2i}) = d^{2g+1} \prod_{j =1}^g\prod_{p \mid d_j}(1-p^{-2j})\, .
    $$
\end{prop}

\subsection{The degree of $\pi_\delta$}\label{section: 72}

As an application, the degree of $\pi_\delta$ (or equivalently, the cardinality of $\operatorname{Sp}(K(\delta))$) is computed:
\begin{prop}\label{prop: degree piD}
    If $\delta = (d_1, \ldots , d_g)$ is a polarization type then
    $$
    \deg (\pi_\delta) = \deg(\phi_\delta) \cdot d_1^{2g-2}d_2^{2g-6}\ldots d_g^{-2g+2}\prod_{1\leq i<j\leq g} \prod_{p \mid d_j/d_i}\frac{(1-p^{-2(j-i)})}{(1-p^{-2(j-i+1)})}\, ,
    $$
    where the product is over primes $p$.
\end{prop}
\begin{proof}
    We will compute this degree by means of the following formulas:
    \begin{itemize}
        \item If $\delta = (1, d_2, \ldots , d_g)$ and $\delta_1 = (d_2, \ldots , d_g)$ then $\deg(\pi_\delta)=\deg (\pi_{\delta_1})$.
        \item If $\delta = (d_1, \ldots , d_g)$ and $\delta_2=(1, d_2/d_1, \ldots , d_g/d_1)$ then
        $$
        \frac{\deg(\pi_\delta)}{\deg (\pi_{\delta_2})} = \frac{\deg(\phi_\delta)}{\deg (\phi_{\delta_2})}\, .
        $$
    \end{itemize}
    The first one follows because the fiber of $\pi_{g,\delta}$ over $(X, \theta)$ is the set of symplectic basis for $\ker (\theta) \cong (\Z^g/\delta\Z^g)^2$, and clearly $(\Z^g/\delta\Z^g)^2 \cong (\Z^g/\delta_1\Z^g)^2$. The second one follows from the Hirzebruch-Mumford Proportionality theorem and the isomorphism
    $$
    \Ac_{g, \delta_2} \cong \Ac_{g,d_1 \cdot \delta_2} = \Ac_{g,\delta}\, .
    $$
    Using these formulas recursively, and noting that when $\delta = (1,\ldots , 1)$ the degrees of $\phi$ and $\pi$ are both $1$, we arrive at the following:
    $$
    \deg (\pi_\delta) = \deg (\phi_{\delta})\cdot \prod_{i=1}^{g-1} \frac{\deg (\phi_{(d_{i+1}/d_i, \ldots , d_g/d_i})}{\deg (\phi_{(1, d_{i+1}/d_i, \ldots , d_g/d_i)})}\, .
    $$
    By Proposition \ref{prop: degree aD},
    $$
    \frac{\deg (\phi_{(d_{i+1}/d_i, \ldots , d_g/d_i})}{\deg (\phi_{(1, d_{i+1}/d_i, \ldots , d_g/d_i)})} = \left(\frac{d_i^{g-i}}{d_{i+1}\cdot \ldots \cdot d_g}\right)^2\cdot \prod_{j=1}^{g-i}\prod_{p \big| \frac{d_{j+i}}{d_i}}\frac{(1-p^{-2j})}{(1-p^{-2j-2})}\, .
    $$
\end{proof}

\subsection{The volume of $T_n$} We compute the volume of $T_n$ using ideas of \cite{D99, BOPY18}

\begin{prop}\label{prop:volumeTn}
    The volume of the Hecke correspondence $T_n$ is given by
    $$
    \prod_{p \mid n}\prod_{i=1}^g\frac{1-p^{g(v_p(n)+i)}}{1-p^{gi}}\, ,
    $$
    where the profuct is over primes.
\end{prop}
\begin{proof}
    Recall that $T_n$ is defined by a correspondence
    $$
\begin{tikzcd}
	& \Bc_{g,n} \\
	{\Ac_g} && {\Ac_g}
	\arrow["\pi_1"', from=1-2, to=2-1]
	\arrow["\pi_2", from=1-2, to=2-3]
\end{tikzcd}
$$
    where both maps are \'etale. In particular, $\operatorname{vol}(T_n)$ equals $\deg(\pi_1)$. For a principally polarized abelian variety $(X,\theta_X)$, the preimages $\pi_1^{-1}[(X, \theta)]$ correspond to isogenies $(X, \theta_X) \longrightarrow (Y, \theta_Y)$ such that $f^* \theta_Y = n \theta_X$.

    The polarization $\theta_X$ induces a Weil pairing on the $n$-torsion points of $X$, and agrees with the pairing that $n\theta$ induces on $\ker (n\theta) = X[n]$. Writing $ \mathbb Z^g /n \mathbb Z^g =\mathbf Z/n$, and $\widehat{H} = \operatorname{Hom}(H, \mathbb C^*)$, it is known that the Weil pairing can be seen as the canonical map
    $$
    \begin{alignedat}{3}
    e :\, &&(\mathbf Z/n\times \widehat{\mathbf Z/n}) \times (\mathbf Z/n\times \widehat{\mathbf Z/n}) &\longrightarrow &\,&\mathbb C^*\\
    &&(x,\widehat{x}) , (y, \widehat{y}) & \mapsto & \, & \widehat{y}(x)\widehat{x}(y)^{-1}
    \end{alignedat}
    $$
    By \cite[Corollary 2.4.4]{BL04}, $\deg(\pi_1)$ agrees with the number of maximal isotropic subgroups $K \subseteq X[n]$. Lemma \ref{lem:debarre} below, which is due to Debarre \cite[proof of Proposition 2.1]{D99} and is stated in \cite[Lemma 1]{BOPY18}, gives a precise description of such subgroups. 

    Suppose that $H = \mathbb Z/d_1\mathbb Z \times \ldots \times \mathbb Z/d_g \mathbb Z$ where $d_i$ divides $d_{i+1}$. An homomoprhism
    $$
    \mathbb Z^g \longrightarrow H
    $$
    factors through $n \mathbb Z^g$ if and only if $d_g$ divides $n$. This homomorphism is given by a matrix in $M_{g \times g}(\mathbb Z)$, which can be uniquely put in Hermite normal form. This is given uniquely by a matrix $(a_{ij})$, where $1 \leq a_{ij} \leq a_{jj}$ for $i <j$, $a_{ij}=0$ for $i >j$ and $a_{ii} = d_i$. There are
    $$
    d_2 d_3^2 \ldots d_g^{g-1}
    $$
    such matrices.
    
    The number of symmetric homomorphisms $u : H \longrightarrow \widehat{H}$ is counted by symmetric matrices, and the total number is
    $$
    d_1^gd_{2}^{g-1} \ldots d_g\, .
    $$

    Putting all of this together, we see that
    $$
    \deg(\pi_1) = \sum_{\substack{(d_1, \ldots , d_g)\text{ such that}\\
    d_1 \mid d_2 \mid \ldots \mid d_g \mid n}}(d_1 \ldots d_g)^{g}\, .
    $$
    If $n = p^k$, the number agrees with
    $$
    \sum_{l \geq 0} p^{lg}\mathcal P(l,k, g)\, ,
    $$
    where $\mathcal P(l,k,g)$ is the number of partitions of $l$ into at most $g$ parts, each of which is bounded by $k$. It is well-known that
    $$
    \sum_{l \geq 0} \mathcal P(l,k,g) q^l = \prod_{i=1}^g\frac{1-q^{k+i}}{1-q^i}\, .
    $$
    After the substitution $q = p^g$ we obtain the expected result when $n$ is a prime power, and in general it follows by decomposing $n$ into prime powers.
\end{proof}

\begin{lem}\label{lem:debarre}
    Given a quotient $\mathbf Z/n \longrightarrow H$ and a symmetric homomorphism $u : H \longrightarrow \widehat{H}$ then
    $$
    K_u = \{(x, \widehat{x}) \in \mathbf Z/n \times H^* : u (x + H) = \widehat{x} \}
    $$
    is a maximal isotropic subgroup of $\mathbf Z/n \times \widehat{\mathbf Z/n}$, and all maximal isotropic subgroups have this form.
\end{lem}

\newcommand{\etalchar}[1]{$^{#1}$}

\end{document}